\documentclass{article}
\usepackage{amsmath}
\usepackage{amsfonts}
\usepackage{amssymb}
\renewcommand{\bf}{\bfseries}
\renewcommand{\it}{\itshape}
\def\a{\alpha}

\def\be{\begin{equation}}
\def\ee{\end{equation}}

\newtheorem{Theorem}{Theorem}[section]

\def\proof#1. {\par
                      \ifdim\lastskip<15pt
                      \removelastskip\penalty-200
                      \vskip5pt plus3pt minus3pt
                      \fi
                       {\def\a{#1}
                       \ifx\a\empty
                       {\noindent\bf Proof.}
                       \else
                       {\noindent\bf Proof of #1.}
                       \fi}\enspace}

\def\endproof{\hfill\hspace{-6pt}\rule[-4pt]{6pt}{6pt}
\vskip8pt plus3pt minus 3pt}

\allowdisplaybreaks[4]

\usepackage{epsf}
\usepackage{epsfig}
\usepackage{floatflt}
\usepackage{float}

 \title{ \Large Orthogonal sequences constructed from quasi-orthogonal  ultraspherical polynomials}

 \author{ \large Oksana Bihun$^\dag$\footnote{Corresponding author: \textit{obihun@uccs.edu}.}~~and Kathy Driver$^\ddag$\\
 \footnotesize
 $^\dag$Department of Mathematics, University of Colorado, Colorado Springs\\
  \footnotesize
  1420 Austin Bluffs Pkwy, Colorado Springs, CO 80918, USA\\
   \footnotesize
 $^\ddag$Department of Mathematics and Applied Mathematics, University of Cape Town, \\
  \footnotesize Private Bag X3, Rondebosch 7701, South Africa}
 \normalsize
 
 \date{}
 
\begin{document}


\maketitle

\hrulefill

\begin{abstract}

\noindent Let $\displaystyle \{x_{k,n-1}\} _{k=1}^{n-1}$   and   $\displaystyle \{x_{k,n}\} _{k=1}^{n},$  $n \in \mathbb{N}$, be two sets of real, distinct points satisfying the interlacing property $ x_{i,n}<x_{i,n-1}< x_{i+1,n}, \, \, \, i = 1,2,\dots,n-1.$   In   \cite {Wen},  Wendroff  proved that  if  $p_{n-1}(x) = \displaystyle \prod _{k=1}^{n-1} (x-x_{k,n-1})$   and  $p_n(x) = \displaystyle \prod _{k=1}^n (x-x_{k,n})$ , then  $p_{n-1}$ and $p_n$ can be embedded in a non-unique monic orthogonal sequence  $\{p_{n}\} _{n=0}^\infty. $   We investigate a question raised by Mourad Ismail at OPSFA 2015 as to the nature and properties of orthogonal sequences generated by applying Wendroff's Theorem to the interlacing zeros of  $C_{n-1}^{\lambda}(x)$ and $ (x^2-1) C_{n-2}^{\lambda}(x)$,  where  $\{C_{k}^{\lambda}(x)\} _{k=0}^\infty$  is a sequence of monic ultraspherical polynomials and $-3/2 < \lambda < -1/2,$  $\lambda \neq -1.$  We construct an algorithm for generating infinite monic orthogonal sequences $\{D_{k}^{\lambda}(x)\} _{k=0}^\infty$ from the two polynomials $D_n^{\lambda} (x): = (x^2-1) C_{n-2}^{\lambda} (x)$ and  $D_{n-1}^{\lambda} (x): = C_{n-1}^{\lambda} (x)$, which is applicable for each pair of fixed parameters $n,\lambda$ in the ranges $n \in \mathbb{N}, n \geq 5$ and  $\lambda > -3/2$, $\lambda \neq -1,0, (2k-1)/2, k=0,1,\ldots$.    We plot and compare the zeros of  $D_m^{\lambda} (x)$ and  $C_m^{\lambda} (x)$ for several choices of  $m \in \mathbb{N}$ and a range of values of the parameters $\lambda$ and $n$. For $-3/2 < \lambda < -1,$ the curves that the zeros of  $D_m^{\lambda} (x)$ and  $C_m^{\lambda} (x)$ approach are substantially different for large values of $m.$ When $-1 < \lambda  < -1/2,$ the two curves have a similar shape while the curves are almost identical for $\lambda >-1/2.$

\vspace{2mm}



\end{abstract}




\section{\large{Introduction}}

\noindent The monic ultraspherical polynomial $C_{n}^{\lambda}(x)$ is defined by the three term recurrence relation \cite[eqn.(8.18)]{KoeSwa}

\begin{equation}
C_n^{\lambda}(x) = x C_{n-1}^{\lambda}(x)- b_n^{\lambda} C_{n-2}^{\lambda}(x),\, \, \,\lambda \neq 0,-1,\dots ; \, \,  n = 1,2,\dots\ ,
\label{eq1} \end{equation}

\noindent where

\begin{equation} C_{-1}^{\lambda}(x)\equiv 0,\,\, \,  C_0^{\lambda}(x)=1, \, \,\, \,  b_n^{\lambda} = \frac{(n-1)(n-2+2\lambda)}{4(n-2+\lambda)(n-1+\lambda)}, \, \, \,\lambda \neq 0,-1,\dots ; \, \,  n = 1,2,\dots  \label{eq2} \end{equation}

\noindent For each $\lambda > -\frac{1}{2},$ the sequence $\{C_{n}^{\lambda}(x)\} _{n=0}^\infty$ is orthogonal on $(-1,1)$ with respect to the weight function $(1-x^2)^{\lambda -\frac{1}{2}}$ and for each $n \in \mathbb{N}, $ $ n \geq 1,$  the zeros of $C_{n}^{\lambda}(x)$  are real, simple, symmetric, lie in $(-1,1)$ and  the zeros of $C_{n-1}^{\lambda}(x)$ interlace with the zeros of $C_{n}^{\lambda}(x),$  $n \geq 2,$ \, ( see \cite[Theorem 3.3.2]{Sze}) namely,

\begin{equation}
 -1 < x_{1,n} < x_{1,n-1} < \dots < x_{n-1,n} < x_{n-1,n-1} < x_{n,n} < 1. \end{equation}

\noindent where $\displaystyle \{x_{i,n}\} _{i=1}^{n}$  are the zeros of $C_{n}^{\lambda}(x)$ in increasing order.

\medskip

\noindent   As $\lambda$ decreases below $-1/2,$ two (symmetric) zeros of $C_{n}^{\lambda}(x)$ leave the interval  $(-1,1)$ through the endpoints $-1$ and $1$ ( see \cite [p. 296]{DrDu} ) and remain real with absolute value $>1$ for each $n \in \mathbb{N},$ $n \geq 3,$ and  $-\frac32 < \lambda < -\frac12, \lambda \neq -1.$  For $-\frac32 < \lambda < -\frac12, \lambda \neq -1,$  the sequence  $\{C_{n}^{\lambda}(x)\}_{n=0}^\infty$ is quasi-orthogonal of order $2$ with respect to the weight function $(1-x^2)^{\lambda +\frac{1}{2}}$ , (see  \cite [Theorem 6] {BDR} and \cite [p.144]{BDR1}) and, for any $n \in \mathbb{N}, n \geq 4,$   ( see \cite[Theorem 3.1] {DrMu1}),

\begin{equation}
 x_{1,n-1} <  x_{1,n} < -1 < x_{2,n} < x_{2,n-1} < \dots < x_{n-2,n-1} < x_{n-1,n} < 1 < x_{n,n} < x_{n-1,n-1}. \end{equation}

\noindent It follows from  $(4)$ that  the zeros of  $C_{n-1}^{\lambda}(x)$ and $C_{n}^{\lambda}(x)$  are not interlacing for any $n \in \mathbb{N}, n \geq 4$ and   $-\frac{3}{2} < \lambda < -\frac{1}{2}, \lambda \neq -1,$  but we see
from $(3)$ and  $(4)$ that the zeros of  $C_{n}^{\lambda}(x)$ interlace with the zeros of $ (x^2-1) C_{n-1}^{\lambda}(x)$ for each $n \in \mathbb{N}, n \geq 4$ and each $\lambda$ with $ \lambda > -\frac {3}{2}, \lambda \neq -1.$
\medskip

\medskip

\noindent  In 1961,  Wendroff  \cite [p. 554] {Wen} proved that, for any fixed positive integer $n,$ $n \geq 2,$ if $\displaystyle \{x_i\} _{i=1}^{n}$ and $\displaystyle \{y_i\} _{i=1}^{n-1}$ are two sets of real, distinct points satisfying the interlacing property $ x_{1}<y_{1}< x_{2}< y_{2}<\dots <x_{n-1} < y_{n-1} <x_{n},$ there exist infinitely many sequences $\{p_{k}(x)\} _{k=0}^\infty$  of monic orthogonal polynomials with $p_n(x) = \displaystyle \prod _{k=1}^n (x-x_k)$ and $p_{n-1}(x) = \displaystyle \prod _{k=1}^{n-1} (x-y_k).$  His proof is constructive and for a given, fixed $n \in \mathbb{N},$ $n \geq 2,$ each polynomial $p_k(x)$ of degree $k \leq n-2$ is uniquely determined by $p_n(x)$ and $p_{n-1}(x).$ In contrast, the monic polynomials of degree $n+1, n+2, \dots$ in any orthogonal sequence that includes  $p_n(x)$ and $p_{n-1}(x)$  are only constrained by the requirement that any infinite sequence of (monic) orthogonal polynomials satisfies a three term recurrence relation of the form

\begin{equation}
p_{n+k}(x) = (x-a_k) p_{n+k-1}(x)- b_k p_{n+k-2}(x),\, a_k \in \mathbb{R}, b_k >0, \,k = 1,2,\ldots.
\end{equation}

\noindent Since there are infinitely many choices of the coefficients $a_k$ and $b_k$ with $a_k  \in \mathbb{R}$  and $ b_k > 0$  for  $k = 1,2,\dots, $ there are infinitely many distinct monic  orthogonal sequences  $\{p_{k}(x)\} _{k=0}^\infty$  that include  $p_n(x)$ and $p_{n-1}(x).$

\medskip

\noindent Here, we fix $n \in \mathbb{N}, n \geq 5$, fix $\lambda,$ $\lambda > -3/2$, $\lambda \neq -1,0, (2k-1)/2, k=0,1,\ldots$, and define
\begin{equation}
D_{n-1}^{\lambda}(x): = C_{n-1}^{\lambda}(x),\,\, \, \, D_{n}^{\lambda}(x): = (x^2-1)C_{n-2}^{\lambda}(x). \end{equation}

\noindent We investigate the properties of the zeros of polynomials in monic orthogonal sequences  $\{D_m^{\lambda}\}_{m=0}^\infty$  generated by the Wendroff process. It is important to emphasize the dependence on the ``starting value'' of $n \in \mathbb{N}$ when generating each monic orthogonal sequence $\{D_m^{\lambda}\}_{m=0}^\infty$ that includes $D_{n-1}^{\lambda}(x)$ and $D_{n}^{\lambda}(x).$  If $n \in \mathbb{N},$  $n \geq 5$ is large, the number of polynomials (namely, $n-2$) that are uniquely determined in every monic orthogonal sequence that includes $D_{n-1}^{\lambda}(x)$ and $D_{n}^{\lambda}(x)$ is correspondingly large whereas, for example, if $n =5$ we have two degrees of freedom when generating each of the monic polynomials of degree $>5$ and exactly $4$ of the polynomials of lower degree are completely determined.  The restriction $n \geq 5$ arises from the fact that when $-\frac{3}{2} < \lambda < -1,$ the quadratic ultraspherical polynomial $C_{2}^{\lambda}(x)$ has two pure imaginary zeros (see \cite [p. 144]{BDR1}).  If we restrict $\lambda$ to the range $ \lambda > -1,$ the results proved here apply for $n \geq 3.$

\section{\large{Notation}}

\noindent For each $m \in \mathbb{N},$ denote
\begin{equation}
C_m^\lambda(x)=x^m+\sum_{j=1}^m \alpha_{j,m} \,  x^{m-j}=\prod_{j=1}^m(x-x_{j,m}),  \,\, \, \, \, \, c_m=\sum_{j=1}^m x_{j,m}.
\end{equation}

\begin{equation}
D_m^\lambda(x)=x^m+\sum_{j=1}^m \beta_{j,m} \, x^{m-j}=\prod_{j=1}^m(x-y_{j,m}),
\, \, \, \, \,  d_m=\sum_{j=1}^m y_{j,m}.\label{2.2}
\end{equation}

\noindent The zeros $\{x_{j,m}\}_{j=1}^m$ of $C_m^\lambda(x)$  are distinct, real and symmetric with respect to the origin for $m \geq 3$ so that $c_m = 0$ for all $m \in \mathbb{N}, m\geq 3$ while the zeros $\{y_{j,m}\}_{j=1}^m$ of $D_m^\lambda(x)$ are distinct, real and symmetric when $m=n$ or $m=n-1$ so that $d_n = d_{n-1} =0.$   Note that
\begin{equation}
y_{1,n} = x_{1,n-2}, \, \, \,  \, y_{2,n} = -1, \,  \, y_{3,n} = x_{2,n-2},  \dots , y_{n-1,n} =1, \, \, \,  \, y_{n,n} = x_{n-2, n-2}
\end{equation}
and
\begin{equation}
y_{k,n-1} = x_{k,n-1},  \, \,  \,\, \, \, \,  k = 1,2,..\dots n-1.
\end{equation}

\section{\large{Orthogonal sequences generated by  $C_{n-1}^{\lambda}(x)$ and $(x^2 -1)C_{n-2}^{\lambda}(x),  \, n \geq 5$}}

\noindent Our main result is the following theorem.

\begin{Theorem}
\label{MainThm}
Let $\{C_m^{\lambda}(x)\}_{m=0}^\infty$ be the sequence of monic ultraspherical polynomials defined by $(1)$ and $(2).$ Fix $n \in{\mathbb N},$ $n \geq 5,$  fix $\lambda \in (-\frac32, -\frac12),$ $\lambda \neq -1,$ and suppose that $\epsilon >0$ is abitrary.

\noindent Define

\begin{equation} a_n:=x_{n-2,n-2} +\epsilon\end{equation}

\noindent  where  $x_{n-2,n-2}>1$ is the largest zero of $C_{n-2}^{\lambda}(x).$
	
\medskip
	
\noindent Let the sequence of  monic polynomials $\{D_m^{\lambda}(x)\}_{m=0}^\infty$ be defined by:

\begin{eqnarray}
D_n^\lambda(x)=(x^2-1)C_{n-2}^\lambda(x), \;\;\;D_{n-1}^\lambda (x)=C_{n-1}^\lambda(x),\\
D_{n-j}^\lambda(x)=-\frac{1}{\ell_{n-j+2}} \left[ D_{n-j+2}^\lambda(x)-xD_{n-j+1}^\lambda(x)\right],\; j=2,3,\ldots,n\\
D_{n+j}^\lambda(x)=xD_{n+j-1}^\lambda(x)-\ell_{n+j}D_{n+j-2}^\lambda(x), \;j=1,2,\ldots \label{Dnpj}
\end{eqnarray}

\noindent where
\begin{eqnarray}
\ell_{n-j}=\beta_{2,n-j-1}-\beta_{2,n-j}>0, \;j=0,1,\ldots,n-2, \\
\ell_{n+j} \in \left(0, \frac{a_n \,D_{n+j-1}^\lambda(a_n)}{D_{n+j-2}^\lambda(a_n)}\right), \;j=1,2,\ldots,
\end{eqnarray}

\noindent and $\beta_{2,m}$ is the coefficient of $x^{m-2}$ in $D_m^\lambda(x)$, see \eqref{2.2}.
Then the sequence $\{D_m^{\lambda}(x)\}_{m=0}^\infty$ is symmetric and orthogonal with respect to a positive measure supported on the interval $(-a_n,a_n)$.

\end{Theorem}

\medskip

\noindent  {\bf{{Proof of Theorem 3.1}}}

\medskip

\noindent Fix $n \in \mathbb{N}, n \geq 5.$  The  monic polynomial $D_{n-2}^\lambda(x)$ is uniquely determined by
\begin{equation} D_n^\lambda(x)-x D_{n-1}^\lambda(x) = -\ell_n  D_{n-2}^\lambda(x) \end{equation}
where $D_n^\lambda(x)$ and $D_{n-1}^\lambda(x)$ are defined by $(12)$ and the coefficient $\ell_n $ is chosen so that $D_{n-2}^\lambda(x)$ is monic.  The positivity of $\ell_{n}$ follows from the interlacing property of the zeros of $D_n^\lambda(x)$ and $D_{n-1}^\lambda(x).$  In the same way, for each $j=3,4,\ldots,n$, the polynomial $D_{n-j}^\lambda(x)$ is constructed from the polynomials $D_{n-j+1}^\lambda(x)$ and $ D_{n-j+2}^\lambda(x.)$   The process is repeated until we obtain $D_{0}^\lambda(x) = 1.$  The polynomials $D_{n+j}^{\lambda}(x),$ $j = 1,2,\dots$  are constructed recursively using the three-term recurrence relation
\begin{equation}
D_{n+j}(x)=x D_{n+j-1}^\lambda(x)-\ell_{n+j} D_{n+j-2}(x), \;\; j=1,2,\ldots,
\end{equation}
\noindent choosing positive coefficients $\ell_{n+j}$ for  $j = 1,2, \dots.$  This ensures (Favard's Theorem) that the infinite sequence  $\{D_m^{\lambda}(x)\}_{m=0}^\infty$ is orthogonal with respect to a positive measure. Wendroff mentions in his proof \cite [p. 554] {Wen} that the coefficients  $\ell_{n+j}$ can be chosen in such a way that all zeros of  $D_{n+j}^\lambda (x)$ lie in the interval $(-a_n,a_n)$ for each $j \geq 1$ but does not indicate how to choose the coefficients to achieve this outcome. Here, we prove that the choice of $\ell_{n+1}$ given in $(16)$ ensures that all the zeros of $D_{n+1}^\lambda (x)$ lie in the interval $(-a_n,a_n).$ Applying the same argument iteratively, it can be shown that the choice of $\ell_{n+j}$ given in $(16)$ ensures that the zeros of  $D_{n+j}^\lambda (x)$ lie in the interval $(-a_n,a_n)$ for each $j \in \mathbb{N}.$

\medskip

\noindent Suppose $\ell_{n+1} \in \left(0, \frac{a_n \,D_{n}^\lambda(a_n)}{D_{n-1}^\lambda(a_n)}\right).$  We show that the zeros of $D_{n+1}^\lambda$ lie in the interval $(-a_n,a_n)$. From $(4)$ with $n$ replaced by $n-1,$  it follows immediately that the zeros of $D_n^\lambda(x) = (x^2 -1) C_{n-2}^\lambda (x) $ and $D_{n-1}^\lambda (x) = C_{n-1}^\lambda (x) $ lie in the open interval $(-a_n,a_n).$ In addition, $D_n^\lambda(x)$ and $D_{n-1}^\lambda (x)$ are monic polynomials with no zeros greater than $x_{n-2,n-2}$ so $D_{n}^\lambda(a_n) >0 $ and $ D_{n-1}^\lambda(a_n)>0.$  Since  $a_n >0,$ it follows that $ \frac{a_n \,D_{n}^\lambda(a_n)}{D_{n-1}^\lambda(a_n)}>0$. By construction, the zeros of  $D_{n+1}^\lambda$ and $D_{n}^\lambda$ interlace and the zeros of $D_n^\lambda$ lie in the interval $(-a_n,a_n)$  so it follows that
$y_{1,n+1}<y_{1,n}<y_{2,n+1}<y_{2,n}<\cdots<y_{n,n+1}<y_{n,n}<y_{n+1,n+1}$
and $-a_n<y_{1,n}<y_{n,n}<a_n.$ We need to show that the largest zero  $y_{n+1,n+1}$ of $D_{n+1}^\lambda$ satisfies $y_{n+1,n+1}<a_n$.

\medskip

\noindent From $(14)$ with $j = 1,$ we have $D_{n+1}^\lambda(x)=x D_n^\lambda(x)-\ell_{n+1} D_{n-1}^\lambda(x) $ so that
\begin{equation}
D_{n+1}^\lambda(y_{n,n})= -\ell_{n+1} D_{n-1}^\lambda(y_{n,n}) <0 .
\end{equation}

\noindent On the other hand, since we are assuming that $\ell_{n+1} \in \left(0, \frac{a_n \,D_{n}^\lambda(a_n)}{D_{n-1}^\lambda(a_n)}\right),$  we have

\begin{equation}
D_{n+1}^\lambda(a_n)=a_n D_n^\lambda(a_n)-\ell_{n+1} D_{n-1}^\lambda (a_n)\, > \, \, a_n D_n^\lambda(a_n)-\frac{a_n D_n^\lambda(a_n)}{D_{n-1}^\lambda(a_n)} D_{n-1}^\lambda (a_n)=0.
\end{equation}

\noindent Therefore,  $D_{n+1}^\lambda(y_{n,n})<0$ and $D_{n+1}^\lambda(a_n)>0$ so that  $y_{n+1,n+1}<a_n$ as required.

\endproof

\medskip

\noindent {\bf{Remark 3.1.}}
The sequences of  polynomials $\{D_m^{\lambda}(x)\}_{m=0}^\infty$ defined in Theorem 3.1 with  $D_{n-1}^\lambda (x) = C_{n-1}^\lambda (x) $  and $D_n^\lambda(x) = (x^2 -1) C_{n-2}^\lambda (x) $  are  (by construction) orthogonal  for each $n\geq 5$ and $\lambda>-\frac{1}{2}$ satisfying $\lambda \neq 0, \lambda \neq \frac{2k-1}{2},$ $k=1,2\ldots$.  It is therefore natural to compare orthogonal sequences  $\{D_m^{\lambda}(x)\}_{m=0}^\infty$ generated by the Wendroff process  with the sequences $\{C_m^\lambda(x)\}_{m=0}^\infty$  of  ultraspherical polynomials orthogonal on $(-1,1)$.   In this case, we can choose  $a_n=1$  so the interval  $(-1,1)$ contains the zeros of all the polynomials in the sequence $\{D_m^\lambda(x)\}_{m=0}^\infty$   as well as  the sequence $\{C_m^\lambda(x)\}_{m=0}^\infty$.

\medskip

\noindent {\bf{Remark 3.2.}}
For $-3/2 < \lambda < -1/2,$ $\lambda \neq -1,$ and $n \geq 5 $ fixed, the largest (real)  zero $x_{n-2,n-2}$ of  $C_{n-2}^\lambda (x) $ is  bounded above by  $  \left( \frac{n-3}{2\lambda +n-2}\right)^{1/2},$  see \cite [(4)]{DrMu2}. An alternative upper bound for  $x_{n-2,n-2}$  is given by\\$1-\frac{2\lambda +1}{(n-2)(n + 2\lambda -2)}, $ see  \cite [(15)]{DrMu2}.
These bounds give estimates for the interval of orthogonality $(-a_n, a_n)$ of the sequences  $\{D_m^\lambda(x)\}_{m=0}^\infty$, where  $D_{n-1}^\lambda (x) = C_{n-1}^\lambda (x) $  and $D_n^\lambda(x) = (x^2 -1) C_{n-2}^\lambda (x) $.

\medskip


\medskip

\noindent {\bf{Remark 3.3.}}
For $-3/2 < \lambda < -1/2,$ $\lambda \neq -1,$ and $n \geq 5 $ fixed, we can choose $\epsilon>0$ in Theorem 3.1  in such a way that the interval $(-a_n, a_n)$ contains all the zeros of the polynomials  $\{C_m^\lambda(x)\}_{m=3}^\infty$.
To this end, we use the estimates $x_{m,m}<\left(\frac{m-1}{2\lambda+m} \right)^{1/2}$ and $x_{m,m}<1-\frac{2\lambda+1}{m(m+2\lambda)}$, $m\geq 3$, see \cite{DrMu2}, where $x_{m,m}$ is the largest zero of $C_m^\lambda$. Because $\max\{ \left(\frac{m-1}{2\lambda+m} \right)^{1/2}: m\geq 3, \lambda > -3/2\}=\left(\frac{2}{2 \lambda+3}\right)^{1/2}$ and $\max\{1-\frac{2\lambda+1}{m(m+2\lambda)}: m\geq 3,\lambda >-3/2\}=\frac{4(2+\lambda)}{3(3+2\lambda)}$, we may choose $a_n := a$ independent of $n \in \mathbb{N}$, namely
\begin{equation}
a=A_1(\lambda):=\left(\frac{2}{2 \lambda+3}\right)^{1/2}
\label{A1}
\end{equation} or
\begin{equation}
a=A_2(\lambda):=\frac{4(2+\lambda)}{3(3+2\lambda)}.
\label{A2}
\end{equation}
To choose the sharper of the two bounds $A_1(\lambda)$ and $A_2(\lambda)$, we use the following comparisons: $A_1(-5/4)=A_2(-5/4)$, $A_1(\lambda)<A_2(\lambda)$ if $\lambda \in (-3/2,-5/4)$ and $A_1(\lambda)>A_2(\lambda)$ if $\lambda \in (-5/4, -1/2)$. Note that $A_1(-1/2)=A_2(-1/2)$ and $\lim\limits_{\lambda \to (-3/2)+} A_k(\lambda)=+\infty$ for $k=1,2$.

\medskip

\noindent {\bf{Remark 3.4.}}
When developing an algorithm for generating orthogonal sequences, we can choose
$\ell_{n+j}=\frac{a D_{n+j-1}^\lambda(a)}{\sigma D_{n+j-2}^\lambda(a)}$, $j=1,2,\ldots$, where $\sigma>1$. Using this  expression for $\ell_{n+j}$ and putting  $x=a$ into \eqref{Dnpj}, we obtain  $\ell_{n+j+1}=\frac{(\sigma-1)}{\sigma^2} a^2$ for all $j=1,2,\ldots$.  The advantage of  choosing $\ell_{n+j}=\frac{a D_{n+j-1}^\lambda(a)}{\sigma D_{n+j-2}^\lambda(a)}$, $j=1,2,\ldots$, with $\sigma>1$, is that all coefficients $\ell_{n+j+1},$  $j=1,2,\ldots$ are equal, namely, $\ell_{n+j+1}=\frac{(\sigma-1)}{\sigma^2} a^2$ for  $j=1,2,\ldots$. This results in a significant reduction in the computational complexity of the algorithm.

\section{\large{Algorithm for construction of  $\{D_m^\lambda(x)\}_{m=0}^{n +k},$  $n$ fixed, $n \geq 5$, $k \in \mathbb{N}$}}

\noindent  Using Theorem~\ref{MainThm} and Remarks 1,3,4, we present an algorithm for construction of the first $n+k+1$ terms of orthogonal sequences $\{D_m^\lambda(x)\}_{m=0}^\infty$.

\begin{enumerate}
\item Choose integers $n\geq 5$ and $k\geq 1$.
\item Choose $\lambda \in (-\frac{3}{2}, +\infty)$ with  $\lambda \neq -1, 0, (2k-1)/2$, $k=0,1,2,\ldots.$
\item \label{choicea} If $-3/2< \lambda < -5/4$, define $a=A_1(\lambda)=\left(\frac{2}{2 \lambda+3}\right)^{1/2}$. If $-5/4 \leq\lambda < -1/2$, define  $a=A_2(\lambda)=\frac{4(2+\lambda)}{3(3+2\lambda)}$. If $\lambda> -\frac{1}{2}$, define $a=1$.
\item Choose  $\sigma>1$.
\item Let $D_n^\lambda(x)=(x^2-1)C_{n-2}^\lambda(x), \;\;\;D_{n-1}^\lambda (x)=C_{n-1}^\lambda(x)$.
\item For $j=0,1$, let $\beta_{2,n-j}$ be the coefficient  of $x^{n-j-2}$ in $D_{n-j}^\lambda(x)$.
\item For $j=2,3,\ldots, n$,\\
let $\ell_{n-j+2}=\beta_{2,n-j+1}-\beta_{2,n-j+2}$; \\ let $D_{n-j}^\lambda(x)=-\frac{1}{\ell_{n-j+2}} \left[ D_{n-j+2}^\lambda(x)-xD_{n-j+1}^\lambda(x)\right]$;\\ let $\beta_{2,n-j}$ be the coefficient  of $x^{n-j-2}$ in $D_{n-j}^\lambda(x)$.
 \item  Let $\ell_{n+1} = \frac{a \,D_{n}^\lambda(a)}{\sigma D_{n-1}^\lambda(a)}$ and
$D_{n+1}^\lambda(x)=xD_{n}^\lambda(x)-\ell_{n+1}D_{n-1}^\lambda(x)$.
\item For $j=2, \ldots, k$,
 let $D_{n+j}^\lambda(x)=xD_{n+j-1}^\lambda(x)-\frac{(\sigma-1)}{\sigma^2} a^2D_{n+j-2}^\lambda(x)$.
 \end{enumerate}
\medskip


\noindent The above algorithm generates the first $n+k+1$ terms of a sequence of symmetric polynomials $\{D_m^\lambda(x)\}_{m=0}^\infty$ orthogonal with respect to some positive measure supported on the interval $(-a,a)$, which contains all the zeros of the symmetric polynomials $\{C_m^\lambda(x)\}_{m=0}^\infty$. If $\lambda \in (-3/2, -1/2)$ and $\lambda \neq -1$, the sequence  $\{C_m^\lambda(x)\}_{m=0}^\infty$ is quasi-orthogonal of order $2$ on $(-1,1)$ with respect to the weight function $(1-x^2)^{\lambda+1/2}$; if $\lambda \in (-1/2, +\infty)$ and $\lambda \neq 0$, $\lambda \neq \frac{2k-1}{2}$ for $k=1,2,\ldots$, the  sequence $\{C_m^\lambda(x)\}_{m=0}^\infty$  is orthogonal with respect to the weight function $(1-x^2)^{\lambda-1/2}$ on the interval $(-1,1)$.

\smallskip

\noindent \textbf{Example 4.1.} In this example we present the first 11 terms of the sequence $\{D_m^\lambda\}_{m=0}^\infty$ using our  algorithm with $n=5$, $k=5$, $\sigma=2$, and $a=\frac{4(2+\lambda)}{3(3+2\lambda)}$, where $\lambda \in (-\frac{3}{2}, +\infty)$,  $\lambda \neq {-1, 0}$   and $ \lambda \neq (2k-1)/2$, $k=0,1,2,\ldots$:

\footnotesize
\begin{eqnarray*}
&&D_0^\lambda(x) = 1,\\
&&D_1^\lambda(x) = x,\\
&&D_2^\lambda(x) = x^2-\frac{2 \lambda^2+7 \lambda+9}{2 (2 \lambda^3+7 \lambda^2+9 \lambda+6)},\\
&&D_3^\lambda(x) = x^3-\frac{3 (2 \lambda+5)}{2 (2 \lambda^2+7 \lambda+9)}x,\\
&&D_4^\lambda(x) = x^4-\frac{3 }{\lambda+3}x^2+\frac{3}{4 (\lambda^2+5 \lambda+6)},\\
&&D_5^\lambda(x) = x^5-\frac{(2 \lambda+7) }{2 \lambda+4}x^3+\frac{3 }{2 \lambda+4} x,\\
&&D_6^\lambda(x) = x^6\\
&&+\frac{-26624 \lambda^8-315136 \lambda^7-1452096 \lambda^6-3030464 \lambda^5-1350544 \lambda^4+6634848 \lambda^3+14325052 \lambda^2+11993936 \lambda+3814971) }{18 (\lambda+2) (2 \lambda+3)^2 (512 \lambda^5+2944 \lambda^4+5208 \lambda^3+4 \lambda^2-8638 \lambda-6429)}x^4\\
&&-\frac{(-8192 \lambda^7-55552 \lambda^6-93408 \lambda^5+238480 \lambda^4+1249616 \lambda^3+2167224 \lambda^2+1773274 \lambda+577325) }{6 (\lambda+2) (2 \lambda+3)^2 (512 \lambda^5+2944 \lambda^4+5208 \lambda^3+4 \lambda^2-8638 \lambda-6429)} x^2\\
&&+\frac{4 (2 \lambda+1)^2 (80 \lambda^3+426 \lambda^2+753 \lambda+442)}{3 (2 \lambda+3)^2 (512 \lambda^5+2944 \lambda^4+5208 \lambda^3+4 \lambda^2-8638 \lambda-6429)},\\
&&
D_7^\lambda(x) = x^7\\
&&+\frac{(-10240 \lambda^8-121088 \lambda^7-561408 \lambda^6-1198624 \lambda^5-656672 \lambda^4+2255736 \lambda^3+5154212 \lambda^2+4387984 \lambda+1408809) }{6 (\lambda+2) (2 \lambda+3)^2 (512 \lambda^5+2944 \lambda^4+5208 \lambda^3+4 \lambda^2-8638 \lambda-6429)}x^5\\
&&+\frac{4096 \lambda^8+78848 \lambda^7+458688 \lambda^6+1074016 \lambda^5+295376 \lambda^4-3734688 \lambda^3-8130836 \lambda^2-7213054 \lambda-2452023}{18 (\lambda+2) (2 \lambda+3)^2 (512 \lambda^5+2944 \lambda^4+5208 \lambda^3+4 \lambda^2-8638 \lambda-6429)}x^3\\
&&-\frac{2 (512 \lambda^6+3328 \lambda^5+7048 \lambda^4+828 \lambda^3-19042 \lambda^2-28747 \lambda-13742) }{3 (2 \lambda+3)^2 (512 \lambda^5+2944 \lambda^4+5208 \lambda^3+4 \lambda^2-8638 \lambda-6429)}x,\\
&&
D_8^\lambda(x) = x^8\\
&&+\frac{-34816 \lambda^8-411392 \lambda^7-1916352 \lambda^6-4161280 \lambda^5-2589488 \lambda^4+6899568 \lambda^3+16600220 \lambda^2+14333968 \lambda+4637883}{18 (\lambda+2) (2 \lambda+3)^2 (512 \lambda^5+2944 \lambda^4+5208 \lambda^3+4 \lambda^2-8638 \lambda-6429)}x^6\\
&&+\Big[\frac{253952 \lambda^{10}+4967424 \lambda^9+36636672 \lambda^8+134987136 \lambda^7+240904128 \lambda^6+28003872 \lambda^5-813980016 \lambda^4}{162 (\lambda+2) (2 \lambda+3)^4 (512 \lambda^5+2944 \lambda^4+5208 \lambda^3+4 \lambda^2-8638 \lambda-6429)}\\
&&-\frac{1823644104 \lambda^3+1962244068 \lambda^2+1102018370 \lambda+259653399}{162 (\lambda+2) (2 \lambda+3)^4 (512 \lambda^5+2944 \lambda^4+5208 \lambda^3+4 \lambda^2-8638 \lambda-6429)} \Big]x^4\\
&&-\frac{8 (6656 \lambda^8+61760 \lambda^7+214784 \lambda^6+252224 \lambda^5-437944 \lambda^4-1922704 \lambda^3-2812378 \lambda^2-1984129 \lambda-566938) }{27 (2 \lambda+3)^4 (512 \lambda^5+2944 \lambda^4+5208 \lambda^3+4 \lambda^2-8638 \lambda-6429)}x^2\\
&&-\frac{16 (\lambda+2)^3 (2 \lambda+1)^2 (80 \lambda^2+266 \lambda+221)}{27 (2 \lambda+3)^4 (512 \lambda^5+2944 \lambda^4+5208 \lambda^3+4 \lambda^2-8638 \lambda-6429)},\\
&&
D_9^\lambda(x) = x^9\\
&&+\frac{-38912 \lambda^8-459520 \lambda^7-2148480 \lambda^6-4726688 \lambda^5-3208960 \lambda^4+7031928 \lambda^3+17737804 \lambda^2+15503984 \lambda+5049339 }{18 (\lambda+2) (2 \lambda+3)^2 (512 \lambda^5+2944 \lambda^4+5208 \lambda^3+4 \lambda^2-8638 \lambda-6429)}x^7\\
&&+\Big[\frac{376832 \lambda^{10}+6912000 \lambda^9+49677312 \lambda^8+182130432 \lambda^7+333265728 \lambda^6+89989248 \lambda^5-952585632 \lambda^4}{162 (\lambda+2) (2 \lambda+3)^4 (512 \lambda^5+2944 \lambda^4+5208 \lambda^3+4 \lambda^2-8638 \lambda-6429)}\\
&&-\frac{2231977416 \lambda^3+2437175184 \lambda^2+1380264434 \lambda+327276231}{162 (\lambda+2) (2 \lambda+3)^4 (512 \lambda^5+2944 \lambda^4+5208 \lambda^3+4 \lambda^2-8638 \lambda-6429)} \Big]x^5\\
&&-\frac{2 (4096 \lambda^9+166912 \lambda^8+1357504 \lambda^7+4568800 \lambda^6+5470096 \lambda^5-8399264 \lambda^4-38672660 \lambda^3-57223262 \lambda^2-40687679 \lambda-11707302) }{81 (2 \lambda+3)^4 (512 \lambda^5+2944 \lambda^4+5208 \lambda^3+4 \lambda^2-8638 \lambda-6429)}x^3\\
&&+\frac{8 (\lambda+2)^3 (512 \lambda^5+1664 \lambda^4-328 \lambda^3-8108 \lambda^2-13238 \lambda-7313)}{27 (2 \lambda+3)^4 (512 \lambda^5+2944 \lambda^4+5208 \lambda^3+4 \lambda^2-8638 \lambda-6429)}x,\\
&&D_{10}^\lambda(x) = x^{10}\\
&&+\frac{(-14336 \lambda^8-169216 \lambda^7-793536 \lambda^6-1764032 \lambda^5-1276144 \lambda^4+2388096 \lambda^3+6291796 \lambda^2+5558000 \lambda+1820265)}{6 (\lambda+2) (2 \lambda+3)^2 (512 \lambda^5+2944 \lambda^4+5208 \lambda^3+4 \lambda^2-8638 \lambda-6429)}x^8\\
&&+\Big[\frac{57344 \lambda^10+1012736 \lambda^9+7164672 \lambda^8+26224384 \lambda^7+48985088 \lambda^6+18933696 \lambda^5-120883088 \lambda^4}{18 (\lambda+2) (2 \lambda+3)^4 (512 \lambda^5+2944 \lambda^4+5208 \lambda^3+4 \lambda^2-8638 \lambda-6429)}\\
&&-\frac{296145544 \lambda^3+327852636 \lambda^2+187090450 \lambda+44609151}{18 (\lambda+2) (2 \lambda+3)^4 (512 \lambda^5+2944 \lambda^4+5208 \lambda^3+4 \lambda^2-8638 \lambda-6429)} \Big]x^6\\
&&-\Big[\frac{4 (\lambda+2) (200704 \lambda^10+5561856 \lambda^9+45776256 \lambda^8+174881280 \lambda^7+305832480 \lambda^6-27524064 \lambda^5) }{729 (2 \lambda+3)^6 (512 \lambda^5+2944 \lambda^4+5208 \lambda^3+4 \lambda^2-8638 \lambda-6429)}\\
&&-\frac{1252067256 \lambda^4+2680180104 \lambda^3+2833083030 \lambda^2+1572495406 \lambda+366899565}{729 (2 \lambda+3)^6 (512 \lambda^5+2944 \lambda^4+5208 \lambda^3+4 \lambda^2-8638 \lambda-6429)}
\Big]x^4\\
&&+\frac{8 (\lambda+2)^3 (45056 \lambda^7+308992 \lambda^6+680928 \lambda^5-126736 \lambda^4-3262160 \lambda^3-6273816 \lambda^2-5263402 \lambda-1726229) }{243 (2 \lambda+3)^6 (512 \lambda^5+2944 \lambda^4+5208 \lambda^3+4 \lambda^2-8638 \lambda-6429)}x^2\\
&&+\frac{64 (\lambda+2)^5 (2 \lambda+1)^2 (80 \lambda^2+266 \lambda+221)}{243 (2 \lambda+3)^6 (512 \lambda^5+2944 \lambda^4+5208 \lambda^3+4 \lambda^2-8638 \lambda-6429)}.
\end{eqnarray*}
\normalsize

\section{ \large{The zeros of $D_m^{\lambda}(x)$ and $C_m^{\lambda} (x)$}}

\noindent In this section, we plot and compare the zeros of $D_m^{\lambda}(x)$ constructed using the algorithm in Section 4 with the zeros of $C_m^\lambda(x)$, where $m=3,4,\ldots, n+k$ and $n\geq 5$, $k\geq 1$ are fixed integers.

\medskip

\noindent \textbf{Example 5.1.}
Let $n=5$ and $\sigma=2.$  Choose $k=5$,  $\lambda=-5/4$ and $a=\frac{4(2+\lambda)}{3(3+2\lambda)}=2.$  The  polynomials $D_m^\lambda, 0\leq m\leq 10,$ are listed below with the approximate values of their  zeros in curly brackets $\{ \}$ :
\begin{eqnarray*}
&&D_0^\lambda(x) = 1,   \{\}\\
&&D_1^\lambda(x) = x,  \{0\}, \\
&&D_2^\lambda(x) = x^2-\frac{18}{19},  \{-0.973329, 0.973329\}, \\
&&D_3^\lambda(x) = x^3-\frac{10 }{9}x,  \{-1.05409, 0,  1.05409\}, \\
&&D_4^\lambda(x) = x^4-\frac{12 }{7}x^2+\frac{4}{7},  \{-1.12303, -0.673114, 0.673114, 1.12303\}, \\
&&D_5^\lambda(x) = x^5-3 x^3+2 x,  \{-1.41421, -1,  0,  1,  1.41421\}, \\
&&D_6^\lambda(x) = x^6-\frac{72 }{17}x^4+\frac{70}{17}x^2-\frac{12}{17},  \{-1.7026, -1.05773, -0.466529, 0.466529, 1.05773 , 1.7026 \}, \\
&&D_7^\lambda(x) = x^7-\frac{89}{17}x^5+\frac{121 }{17}x^3-\frac{46 }{17}x,  \{-1.83123, -1.10502, -0.812906, 0,  0.812906 , 1.10502 , 1.83123 \}, \\
&&D_8^\lambda(x) = x^8-\frac{106 }{17}x^6+\frac{193 }{17}x^4-\frac{116 }{17}x^2+\frac{12}{17},
\{-1.89282, -1.23417, -1,  -0.359651, 0.359651, 1,  1.23417 , 1.89282 \}, \\
&&D_9^\lambda(x) = x^9-\frac{123 }{17}x^7+\frac{282 }{17}x^5-\frac{237 }{17}x^3+\frac{58 }{17}x,  \\
&&\{-1.92625, -1.41421, -1.05407, -0.643268, 0,  0.643268 , 1.05407 , 1.41421, 1.92625\}, \\
&&D_{10}^\lambda(x) = x^{10}-\frac{140 }{17}x^8+\frac{388 }{17}x^6-\frac{430 }{17}x^4+\frac{174 }{17}x^2-\frac{12}{17}, \\
&& \{-1.94625, -1.55305, -1.09439, -0.867151, -0.292897, 0.292897, 0.867151, 1.09439, 1.55305, 1.94625\}.
\end{eqnarray*}

\medskip

\noindent Note that the zeros of $D_5^{-5/4}(x) $ are  $- \sqrt{2}; -1, 0, 1,  \sqrt{2}.$ The largest and smallest zeros of $D_{10}^{-5/4}(x)$  are close to the limits $-2$ and $2.$

\medskip

\noindent In Figures~\ref{Ex1Zeros3} through~\ref{Ex1Zeros10}, the $y$-coordinates of the plotted points are the zeros of   $D_m^{-5/4}(x)$ (diamond, brown)  and $C_m^{-5/4}(x)$  (round, blue) for $m=3,4,5,10$. The figures suggest that the greatest difference between the zeros of $D_m^{-5/4}(x)$ and $C_m^{-5/4}(x)$ are at the extreme zeros.

\begin{minipage}{\linewidth}
      \centering
      \begin{minipage}{0.45\linewidth}
          \begin{figure}[H]
              \includegraphics[width=\linewidth]{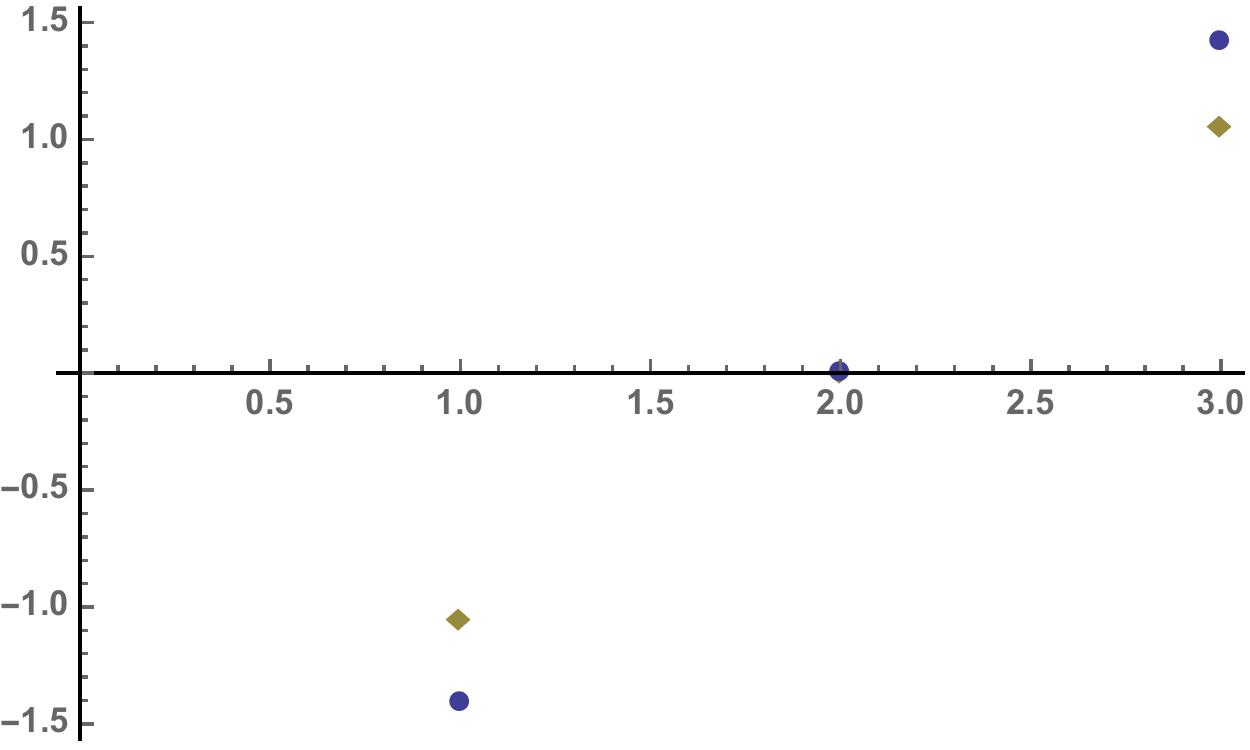}
              \caption{$n=5, \lambda=-5/4, m=3$. The polynomials $C_3^{-5/4}$ and $D_3^{-5/4}$  have a common zero at the origin.}
              \label{Ex1Zeros3}
          \end{figure}
      \end{minipage}
      \hspace{0.05\linewidth}
      \begin{minipage}{0.45\linewidth}
          \begin{figure}[H]
              \includegraphics[width=\linewidth]{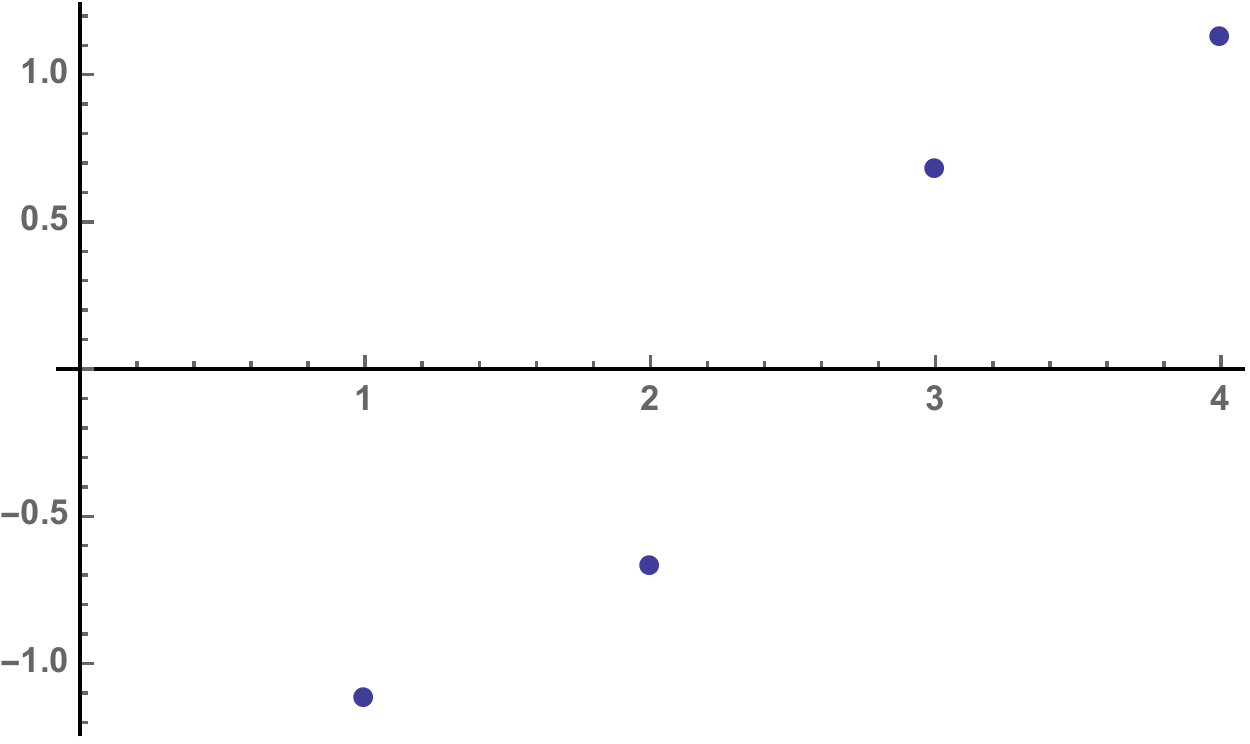}
              \caption{$n=5, \lambda=-5/4, m=4$. Since  $D_4^{-5/4}(x) = C_4^{-5/4}(x)$,  their zeros are equal.  }
              \label{Ex1Zeros4}
          \end{figure}
      \end{minipage}
  \end{minipage}

  \begin{minipage}{\linewidth}
      \centering
      \begin{minipage}{0.45\linewidth}
          \begin{figure}[H]
              \includegraphics[width=\linewidth]{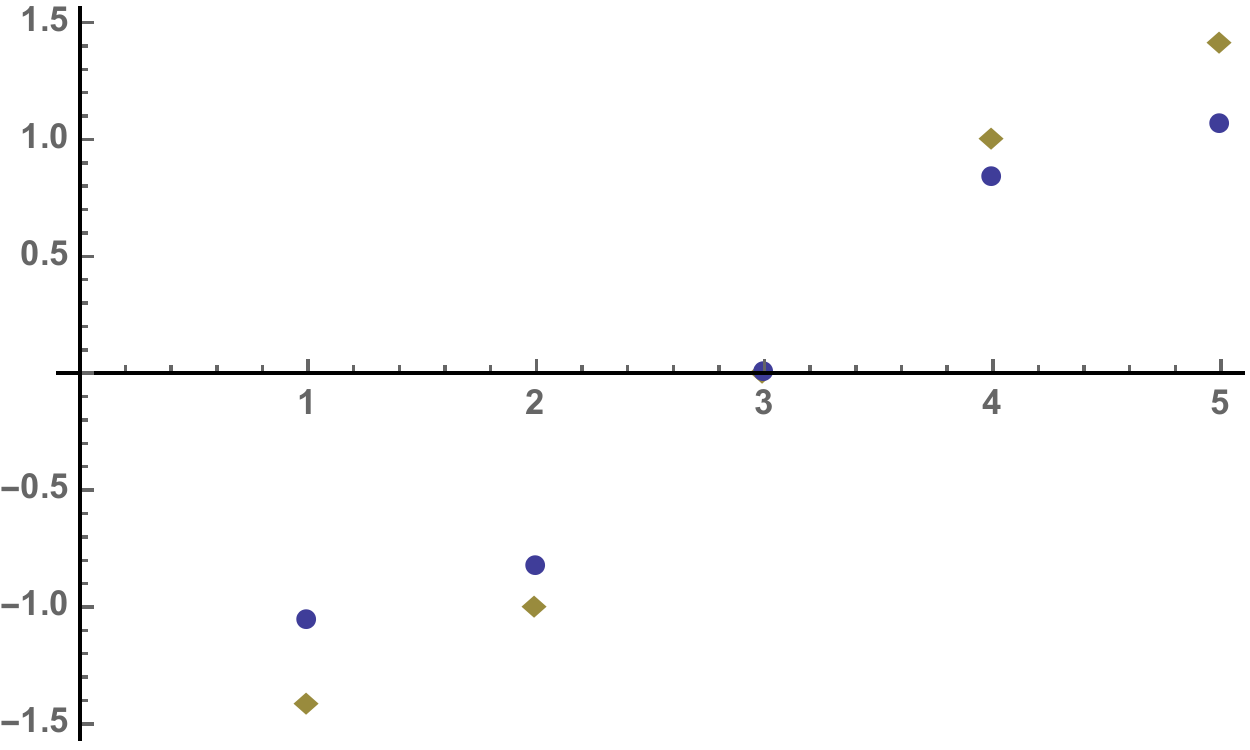}
              \caption{$n=5, \lambda=-5/4, m=5$. By construction, the zeros of $D_5^{-5/4}(x)$ are the zeros of $C_3^{-5/4}(x)$ together with the points $-1$~and~$1$. }
              \label{Ex1Zeros5}
          \end{figure}
      \end{minipage}
      \hspace{0.05\linewidth}
      \begin{minipage}{0.45\linewidth}
          \begin{figure}[H]
              \includegraphics[width=\linewidth]{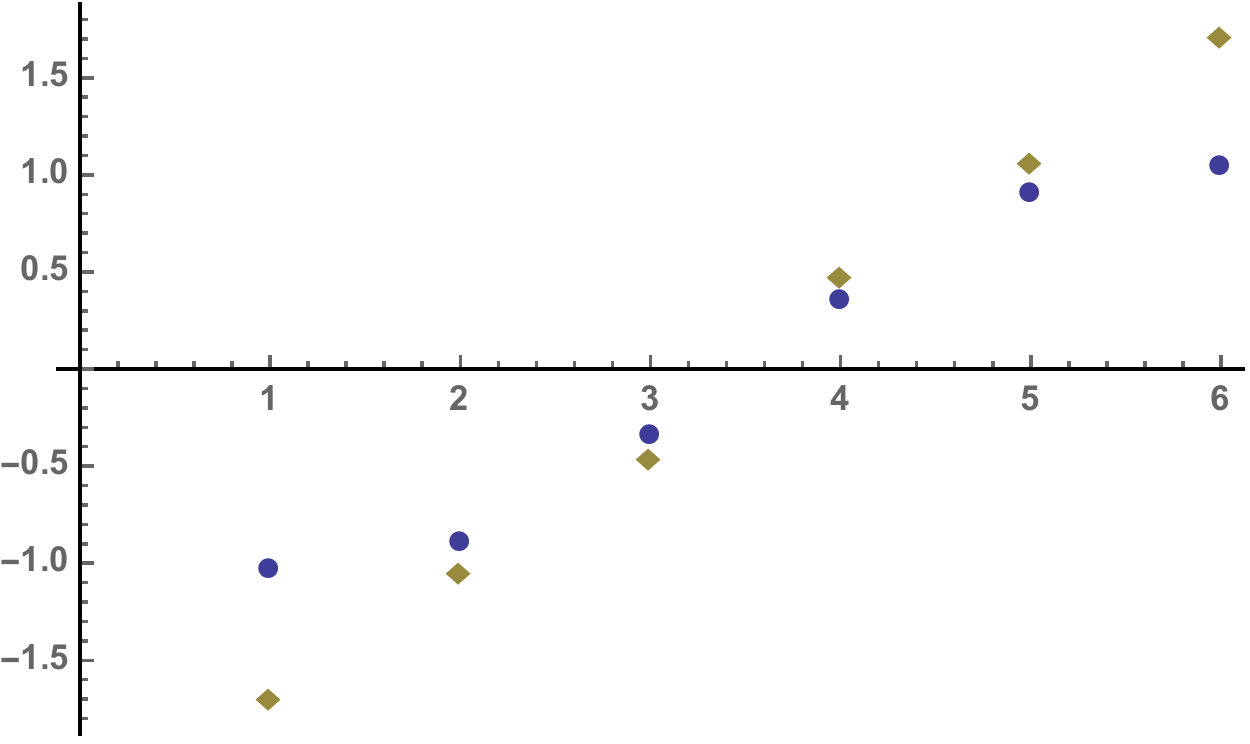}
              \caption{$n=5, \lambda=-5/4, m=10$.}
              \label{Ex1Zeros10}
          \end{figure}
      \end{minipage}
  \end{minipage}

\medskip

\noindent Example 5.1 provides numerical confirmation that the relative ordering of the zeros of $D_{n+1}^{-5/4}, D_{n}^{-5/4},$ and $D_{n-1}^{-5/4}$, is consistent with \cite[Theorem 4]{BDJ}. Replacing  $n$ by $n-1$ and putting $b_n =0$ in (7) and (8) in \cite[Theorem 4]{BDJ}, the negative zeros of $D_{m}^{-5/4},$ $m \in \{n-1,n,n+1\},$ should satisfy
\begin{equation*}
y_{1,n+1} <y_{1,n} <y_{1,n-1} <y_{2,n+1} <y_{2,n} < y_{2,n-1}  \dots
\end{equation*}
while  the  positive zeros of $D_{m}^{-5/4},$ $m \in \{n-1,n,n+1\},$  should  satisfy
\begin{equation*}
y_{n+1,n+1} > y_{n,n}  > y_{n-1,n-1}  > y_{n ,n +1} > y_{n-1 , n }\dots
\end{equation*}
From Example 5.1, we see that the  zeros  of $D_{4}^{-5/4}$, $D_{5}^{-5/4}$,  and $D_{6}^{-5/4}$  satisfy
\begin{equation*}
 y_{1,6}< y_{1,5}< y_{1,4}< y_{2,6} < y_{2,5}< y_{2,4}< y_{3,6}< y_{3,5} =0
 \end{equation*}
 and
 \begin{equation*}
 y_{6,6} > y_{5,5} > y_{4,4} > y_{5,6}> y_{4, 5}> y_{3,4}> y_{4,6}>y_{3,5}=0.
   \end{equation*}
   
\noindent  as expected.

\medskip

\noindent In the examples and figures that follow, we plot the zeros of $D_m^{\lambda}(x)$ and $C_m^{\lambda} (x)$ for selected values of $n$ (the~``starting value'' ), $m$ and $\lambda.$

\smallskip

\noindent \textbf{Example 5.2.}
Let $n=5$, $k=5$, and $\sigma=2$, as in the previous example.
In Figures~\ref{Ex3Zeros10} and~\ref{Ex4Zeros10}, the $y$-coordinates of the plotted points are the zeros of $D_{10}^\lambda$ (diamond, brown) and $C_{10}^\lambda$ (round, blue) respectively for $\lambda=-3/4$, $a=\frac{4(2+\lambda)}{3(3+2\lambda)}=\frac{10}{9}$ and $\lambda=-1/4$, $a=\frac{4(2+\lambda)}{3(3+2\lambda)}=\frac{14}{15}$. 

  \begin{minipage}{\linewidth}
      \centering
      \begin{minipage}{0.45\linewidth}
          \begin{figure}[H]
              \includegraphics[width=\linewidth]{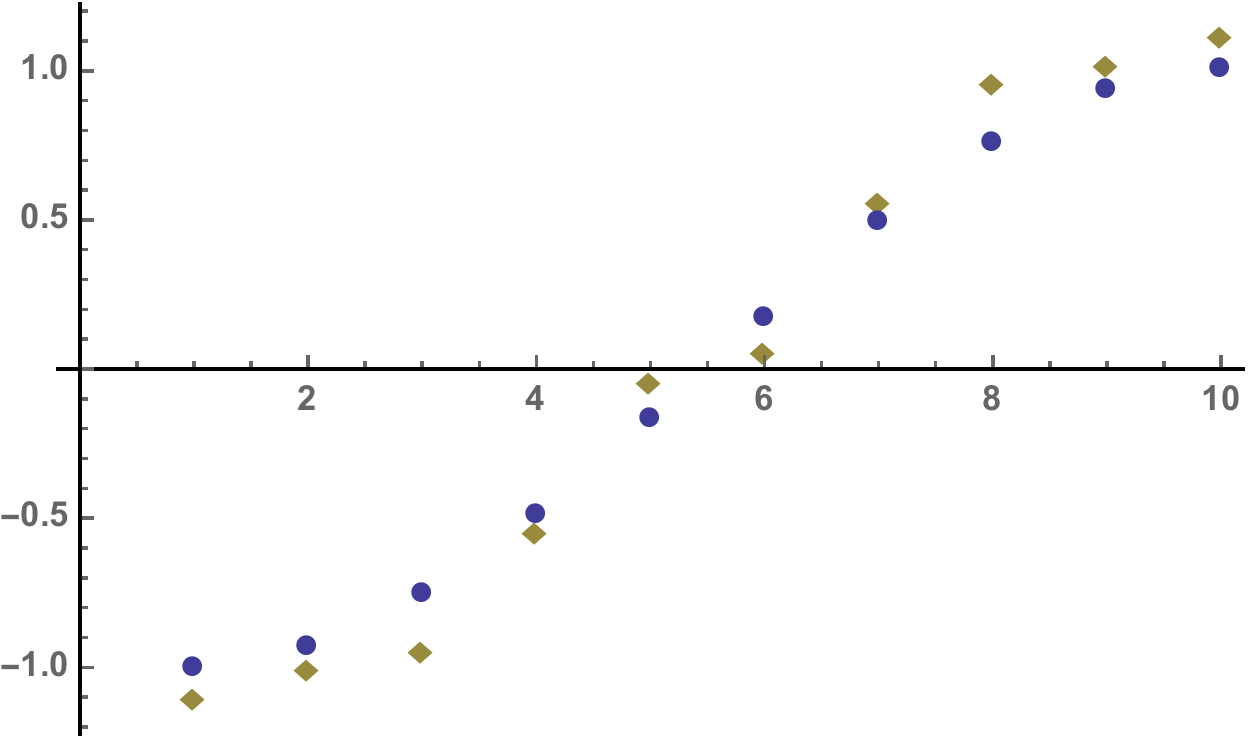}
              \caption{  $n=5, m=10, \lambda=-3/4$. }
              \label{Ex3Zeros10}
          \end{figure}
      \end{minipage}
      \hspace{0.05\linewidth}
      \begin{minipage}{0.45\linewidth}
          \begin{figure}[H]
              \includegraphics[width=\linewidth]{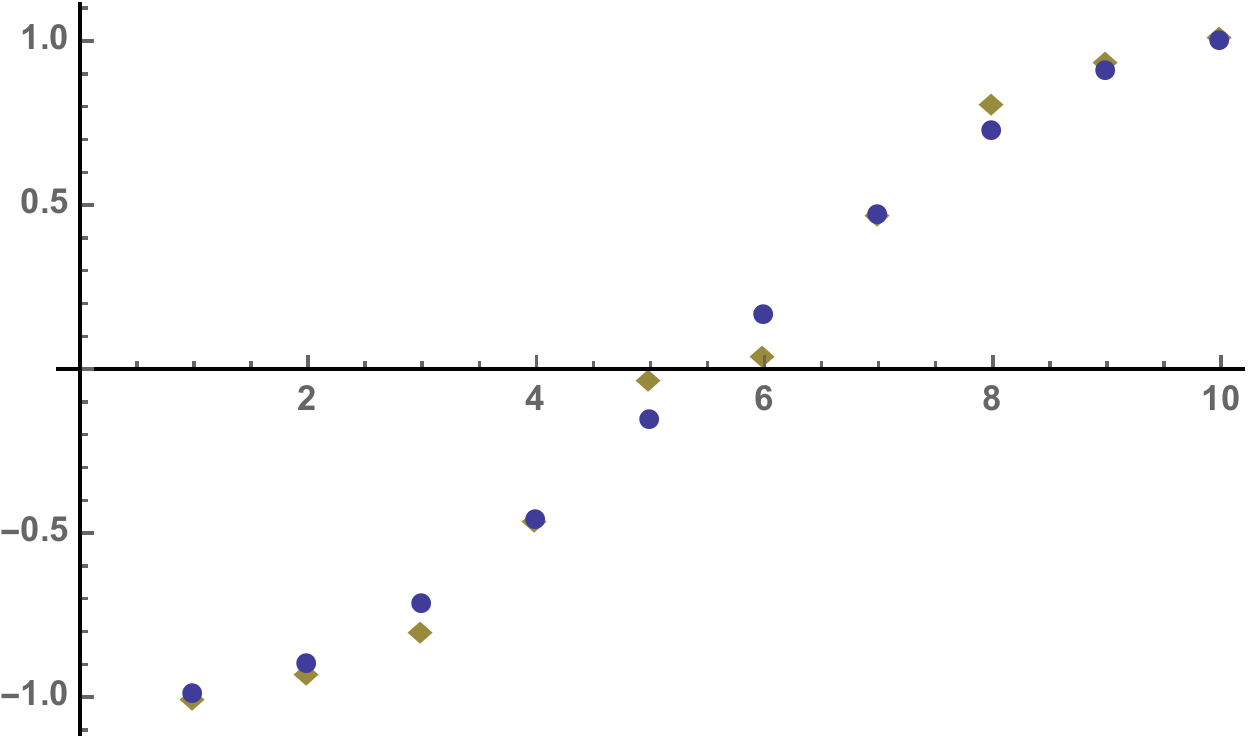}
              \caption{$n=5, m=10, \lambda=-1/4$.  }
              \label{Ex4Zeros10}
          \end{figure}
      \end{minipage}
  \end{minipage}

\smallskip

  \noindent \textbf{Example 5.3.}
Let $n=10$, $k=58$, $\sigma=2.$  Choose $\lambda=-5/4$ and $a=\frac{4(2+\lambda)}{3(3+2\lambda)}=2$.
In Figures~\ref{Ex5Zeros3} through~\ref{Ex5Zeros68} the $y$-coordinates of the plotted points are  the zeros of $D_m^\lambda$ (diamond, brown) and  $C_m^\lambda$ (round, blue) for selected integer values of $m$ between $3$ and $68$. The figures suggest that, as $m$ increases, the curves that fit the zeros of $D_m^{-5/4}$ and $C_m^{-5/4}$ are significantly different.

\begin{minipage}{\linewidth}
      \centering
      \begin{minipage}{0.45\linewidth}
          \begin{figure}[H]
              \includegraphics[width=\linewidth]{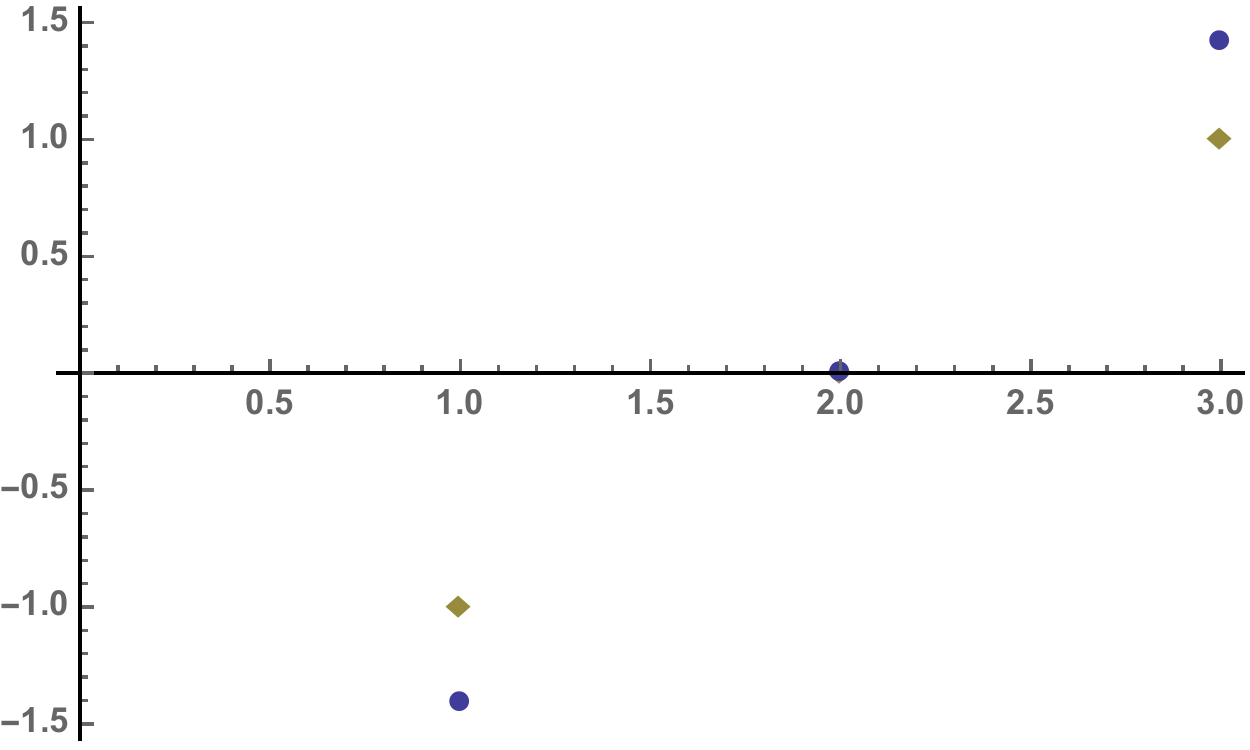}
              \caption{$n=10, \lambda=-5/4, m=3$. }
              \label{Ex5Zeros3}
          \end{figure}
      \end{minipage}
      \hspace{0.05\linewidth}
      \begin{minipage}{0.45\linewidth}
          \begin{figure}[H]
              \includegraphics[width=\linewidth]{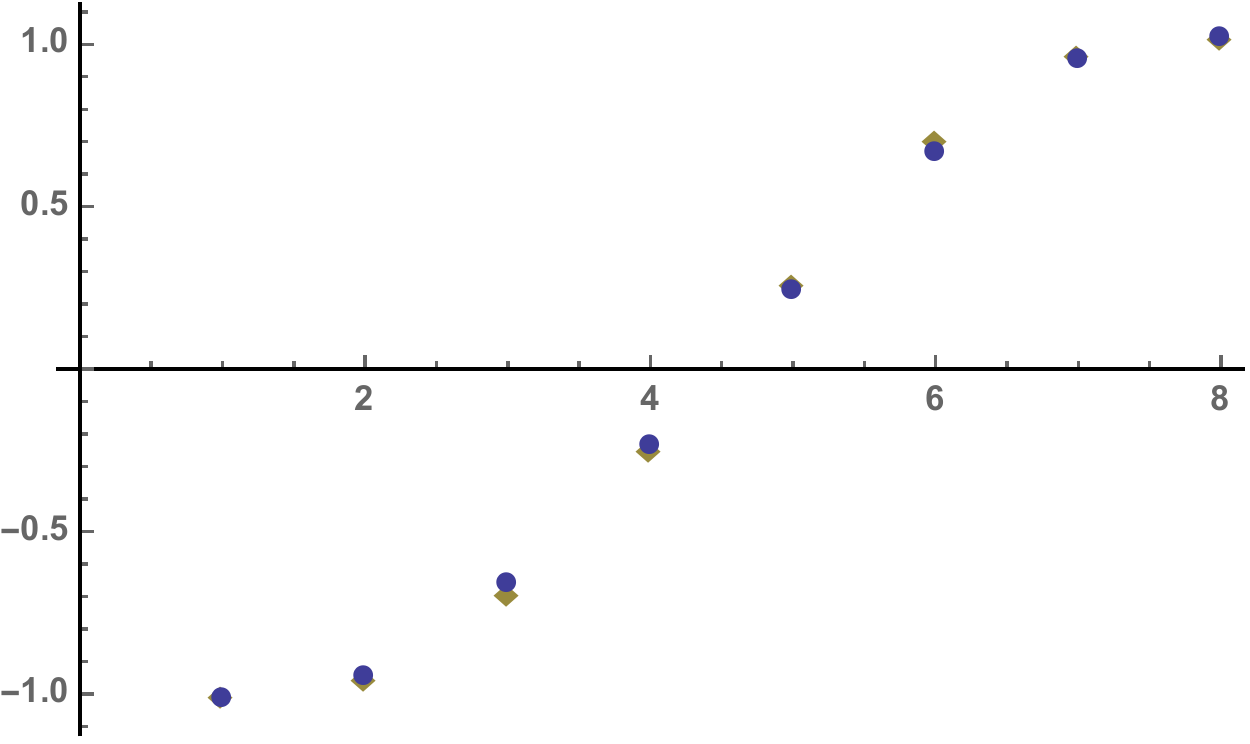}
              \caption{$n=10, \lambda=-5/4, m=8$.  }
              \label{Ex5Zeros8}
          \end{figure}
      \end{minipage}
  \end{minipage}

   \begin{minipage}{\linewidth}
      \centering
      \begin{minipage}{0.45\linewidth}
          \begin{figure}[H]
              \includegraphics[width=\linewidth]{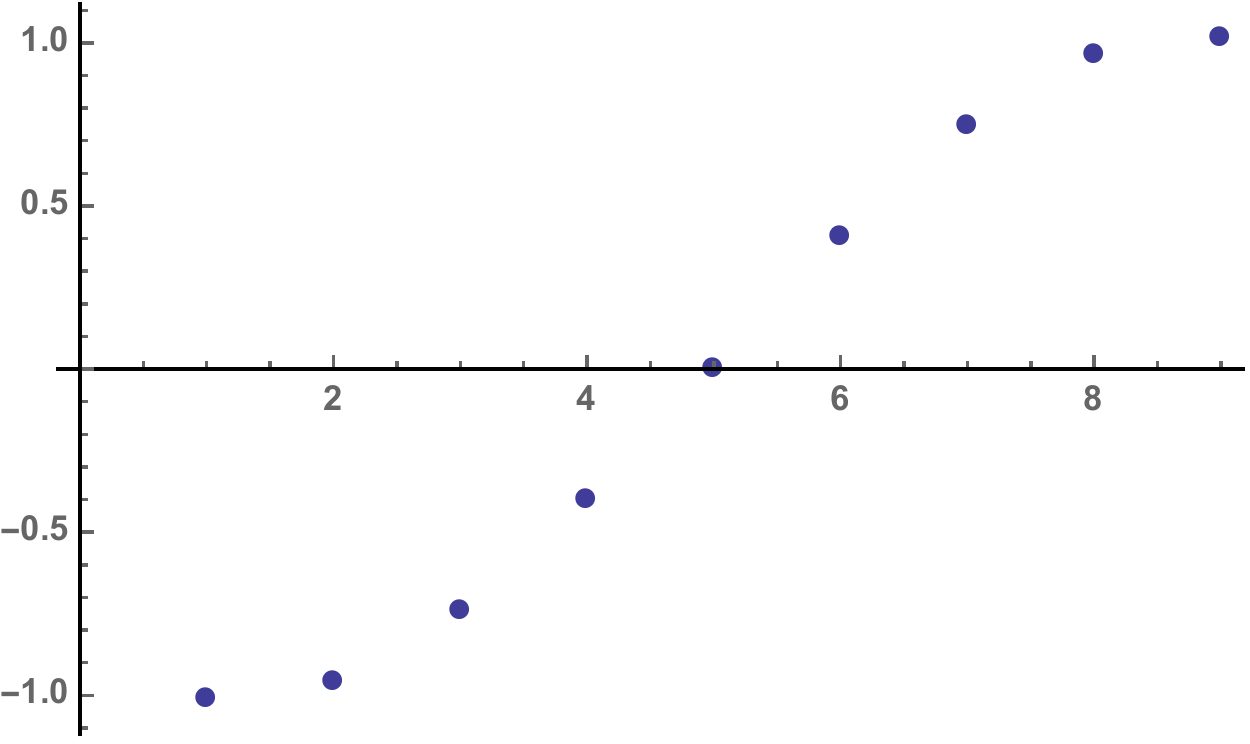}
              \caption{$n=10, \lambda=-5/4, m=9$.}
              \label{Ex5Zeros9}
          \end{figure}
      \end{minipage}
      \hspace{0.05\linewidth}
      \begin{minipage}{0.45\linewidth}
          \begin{figure}[H]
              \includegraphics[width=\linewidth]{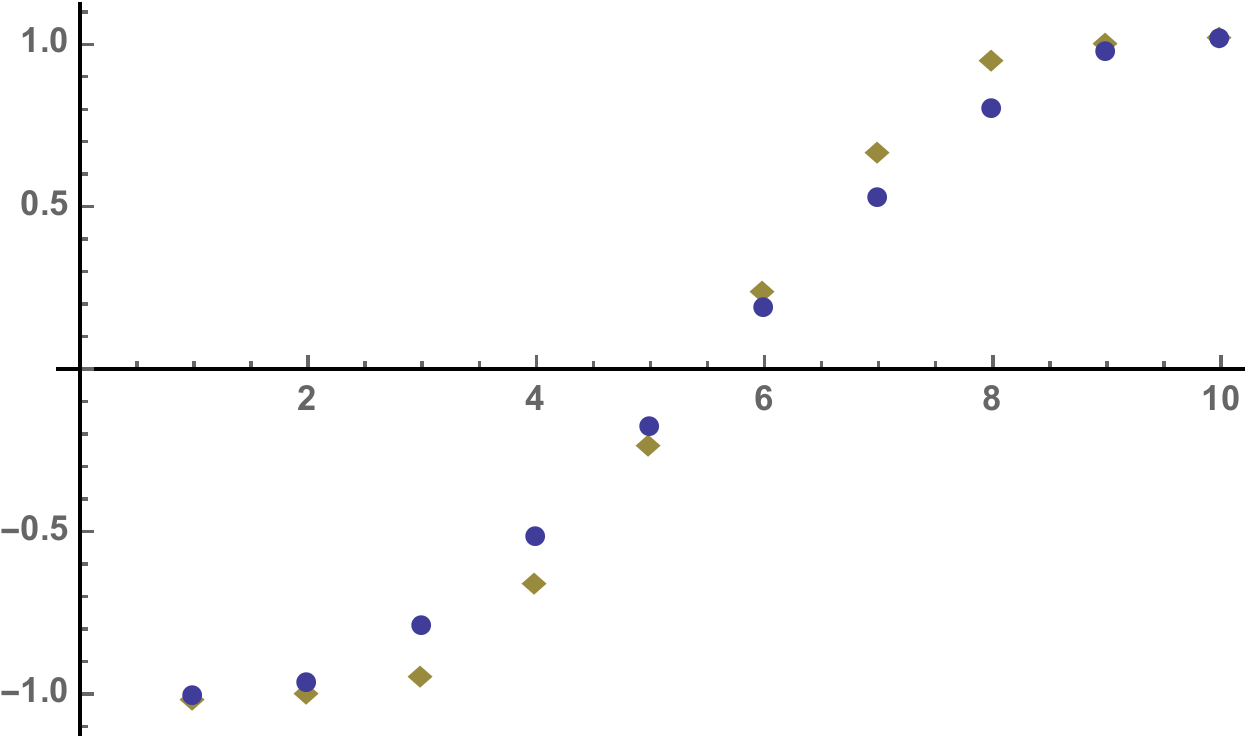}
              \caption{$n=10, \lambda=-5/4, m=10$.}
              \label{Ex5Zeros10}
          \end{figure}
      \end{minipage}
  \end{minipage}

   \begin{minipage}{\linewidth}
      \centering
      \begin{minipage}{0.45\linewidth}
          \begin{figure}[H]
              \includegraphics[width=\linewidth]{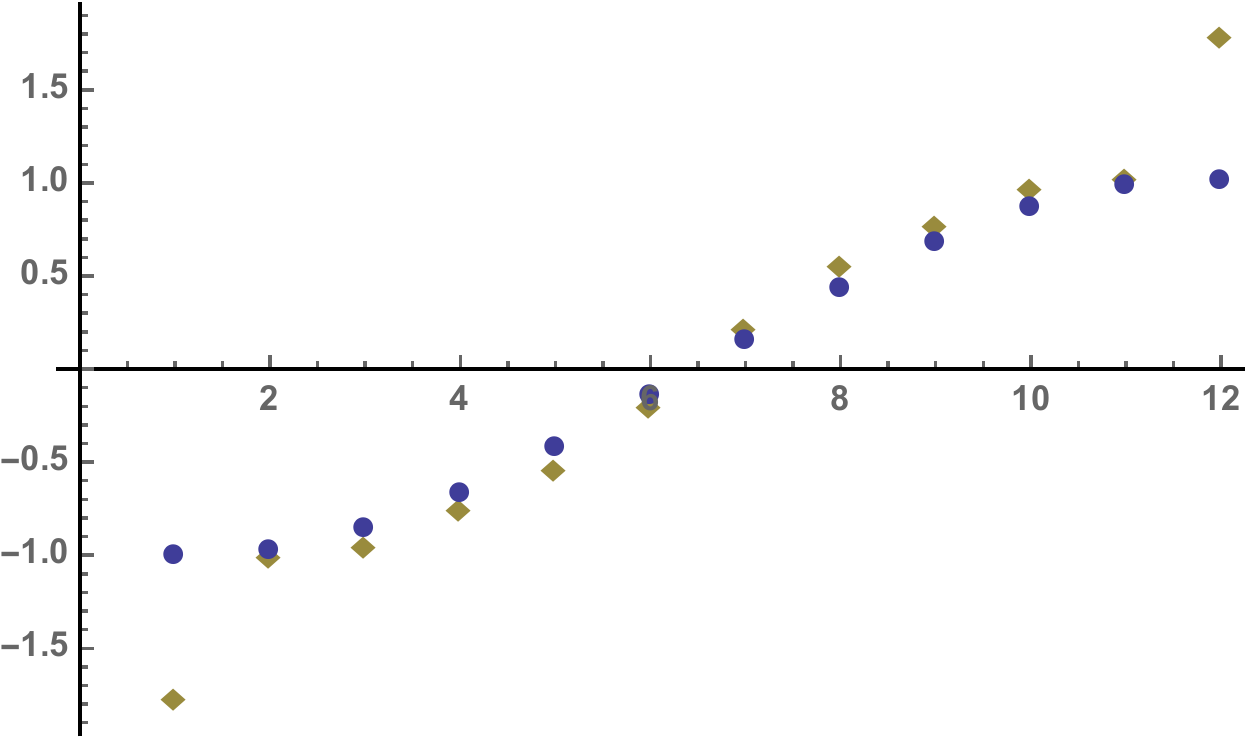}
              \caption{$n=10, \lambda=-5/4, m=12$.}
              \label{Ex5Zeros12}
          \end{figure}
      \end{minipage}
      \hspace{0.05\linewidth}
      \begin{minipage}{0.45\linewidth}
          \begin{figure}[H]
              \includegraphics[width=\linewidth]{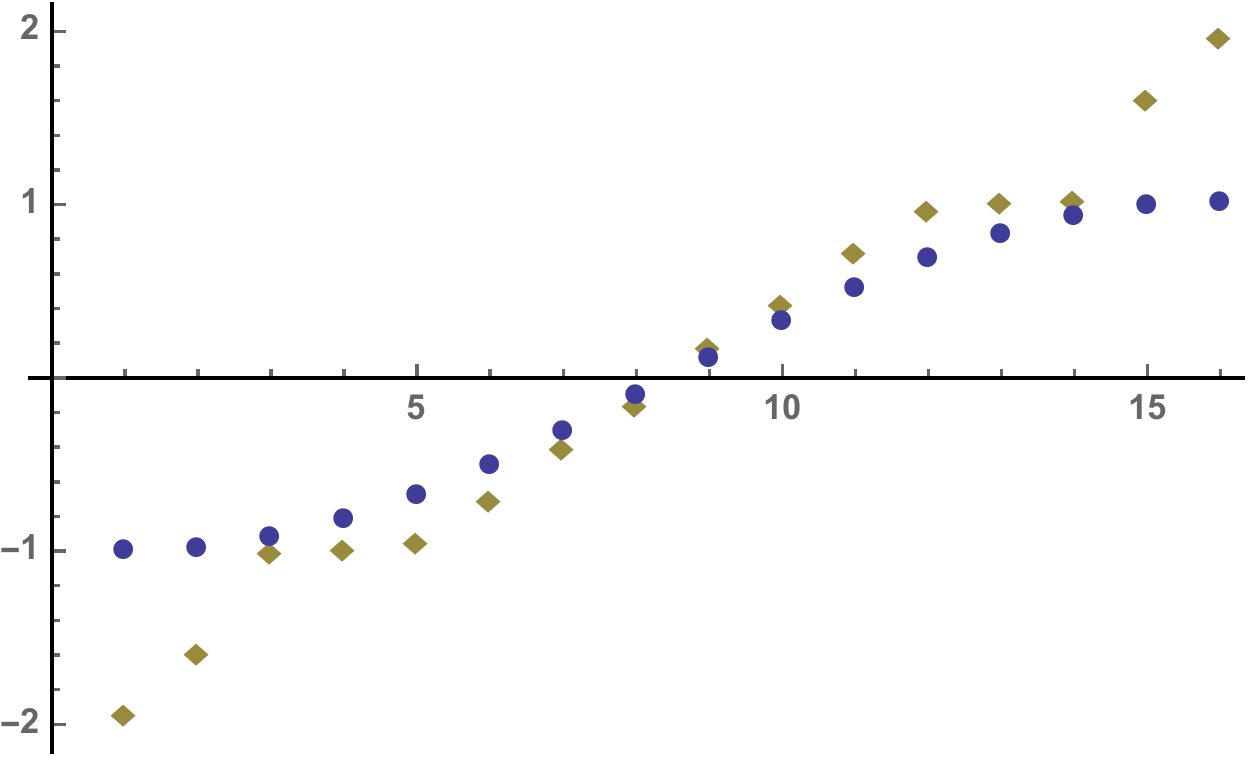}
              \caption{$n=10, \lambda=-5/4, m=16$.  }
              \label{Ex5Zeros16}
          \end{figure}
      \end{minipage}
  \end{minipage}

  \begin{minipage}{\linewidth}
      \centering
      \begin{minipage}{0.45\linewidth}
          \begin{figure}[H]
              \includegraphics[width=\linewidth]{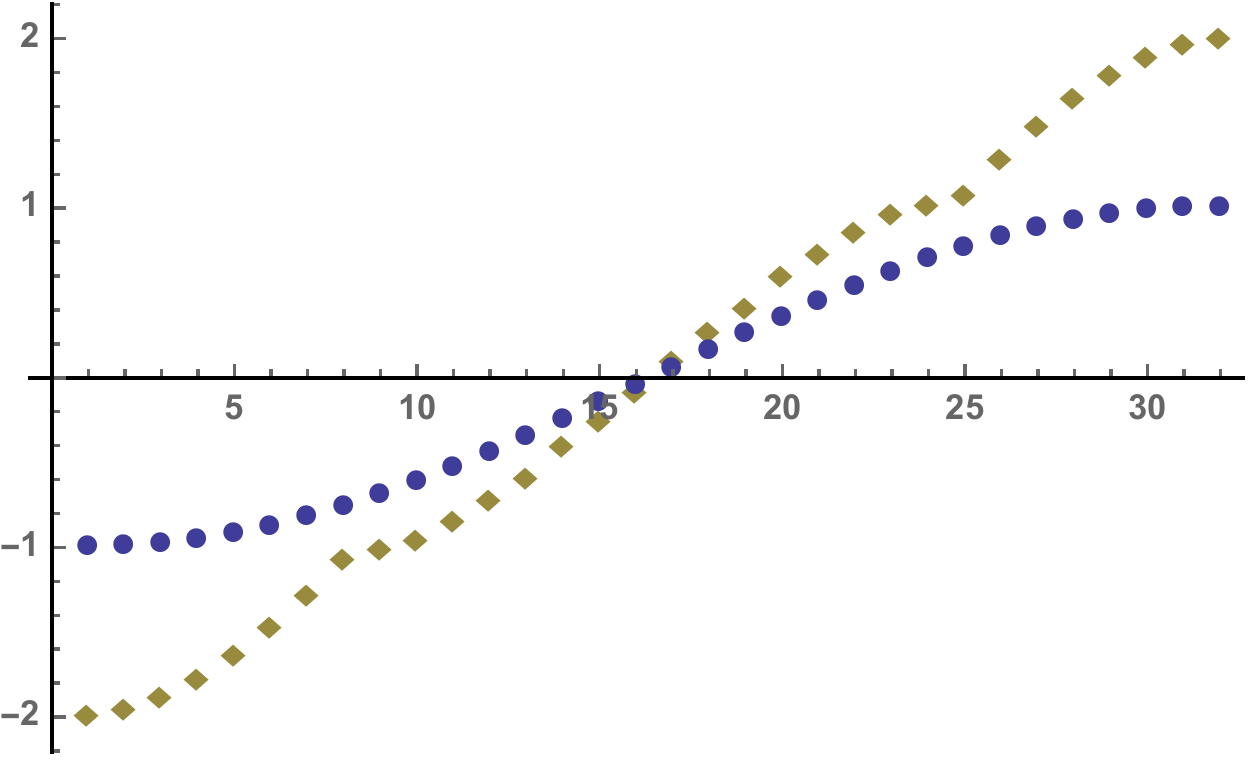}
              \caption{$n=10, \lambda=-5/4, m=32$.}
              \label{Ex5Zeros32}
          \end{figure}
      \end{minipage}
      \hspace{0.05\linewidth}
      \begin{minipage}{0.45\linewidth}
          \begin{figure}[H]
              \includegraphics[width=\linewidth]{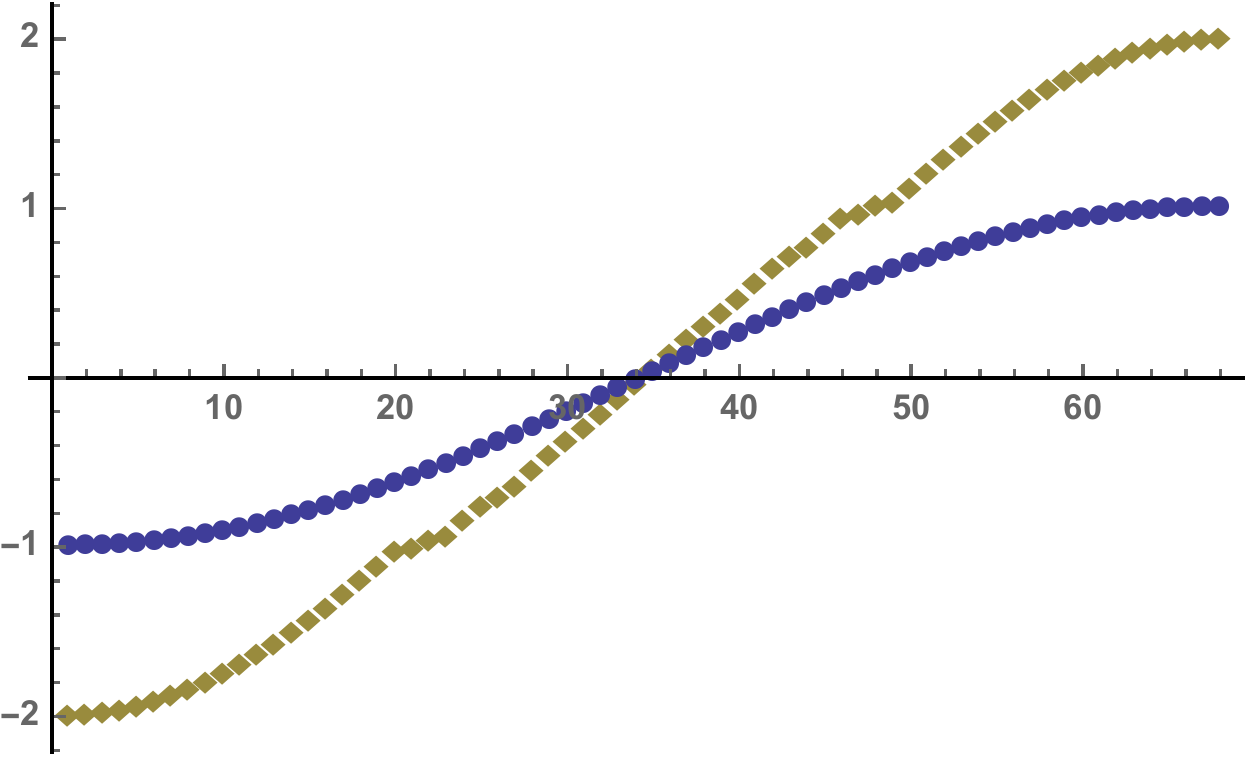}
              \caption{$n=10, \lambda=-5/4, m=68$.}
              \label{Ex5Zeros68}
          \end{figure}
      \end{minipage}
  \end{minipage}

\smallskip

  \noindent \textbf{Example 5.4.}
Let $n=10$, $k=58$, $\sigma=2.$  Choose $\lambda=-3/4$ and $a=\frac{4(2+\lambda)}{3(3+2\lambda)}=\frac{10}{9}$.
In Figures~\ref{Ex6Zeros8} through~\ref{Ex6Zeros68} the $y$-coordinates of the plotted points are the zeros of $D_m^\lambda$ and  $C_m^\lambda$ for selected integer values of $m$ between $8$ and $68$.

    \begin{minipage}{\linewidth}
      \centering
      \begin{minipage}{0.45\linewidth}
          \begin{figure}[H]
              \includegraphics[width=\linewidth]{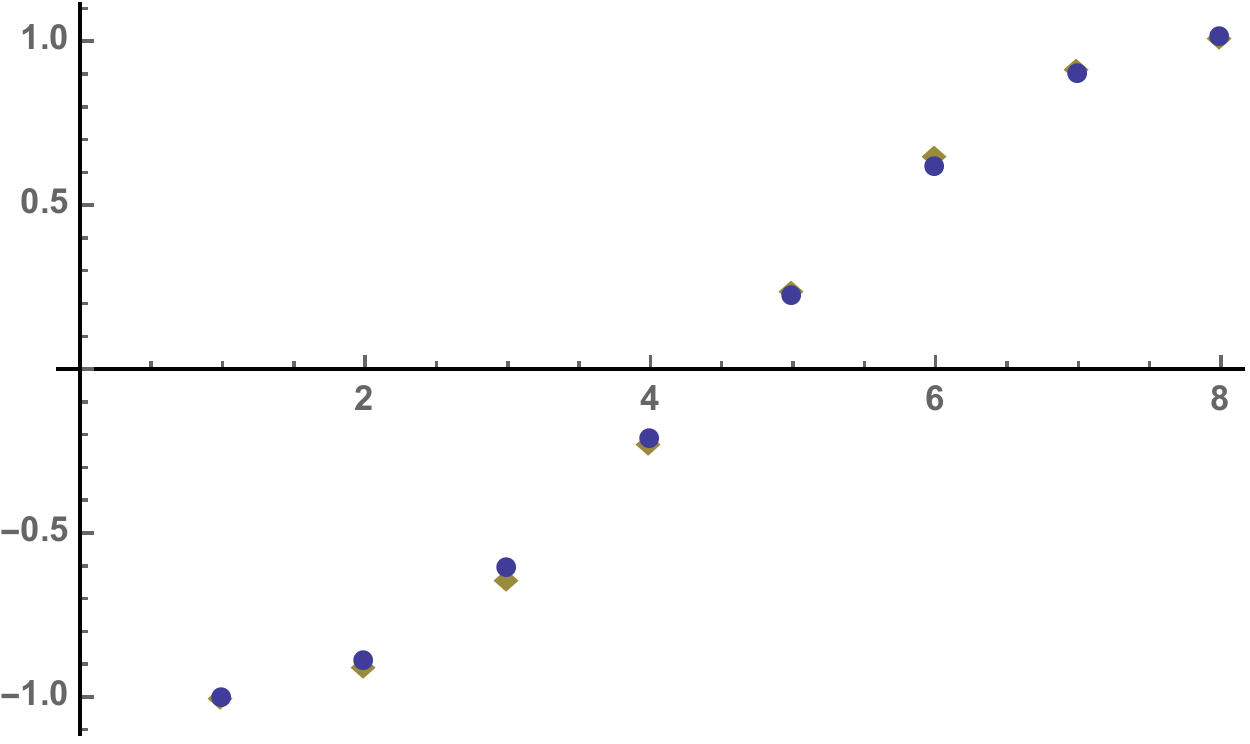}
              \caption{$n=10, \lambda=-3/4, m=8$.}
              \label{Ex6Zeros8}
          \end{figure}
      \end{minipage}
      \hspace{0.05\linewidth}
      \begin{minipage}{0.45\linewidth}
          \begin{figure}[H]
              \includegraphics[width=\linewidth]{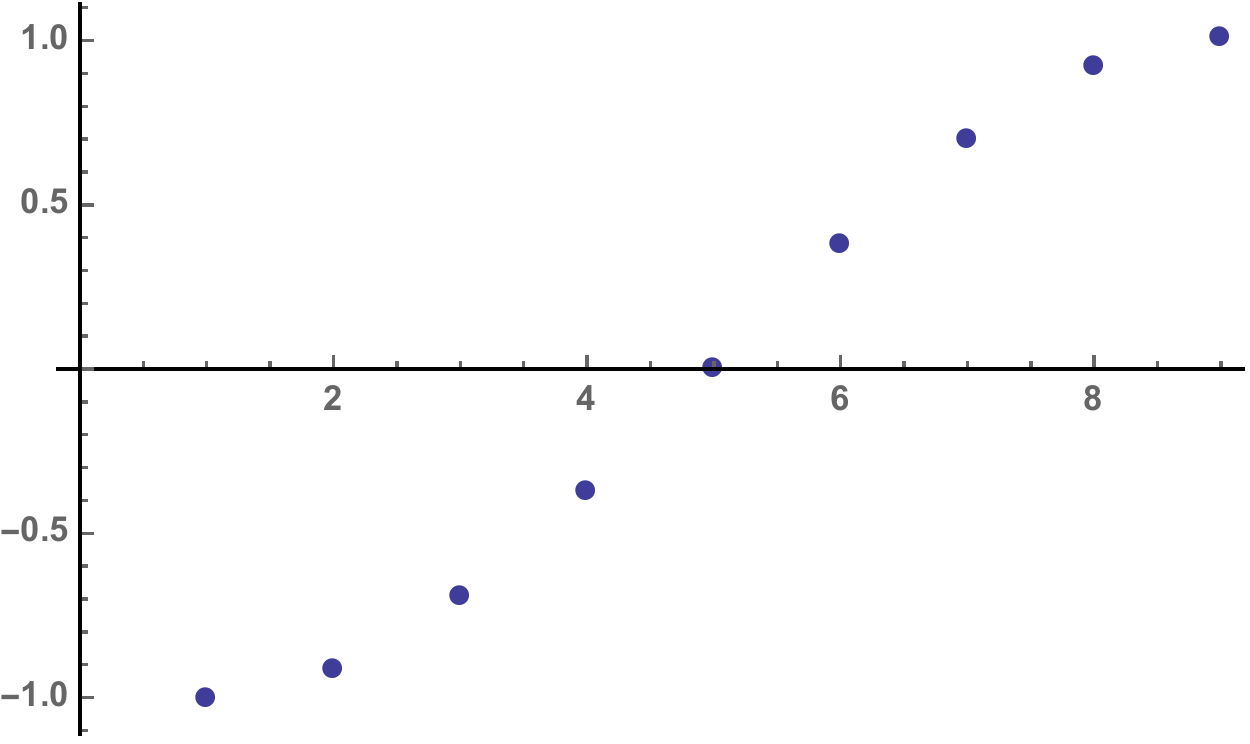}
              \caption{$n=10, \lambda=-3/4, m=9$. Since $C_{9}^{-3/4}(x)=D_{9}^{-3/4}(x)$,  their zeros are equal. }
              \label{Ex6Zeros9}
          \end{figure}
      \end{minipage}
  \end{minipage}

  \begin{minipage}{\linewidth}
      \centering
      \begin{minipage}{0.45\linewidth}
          \begin{figure}[H]
              \includegraphics[width=\linewidth]{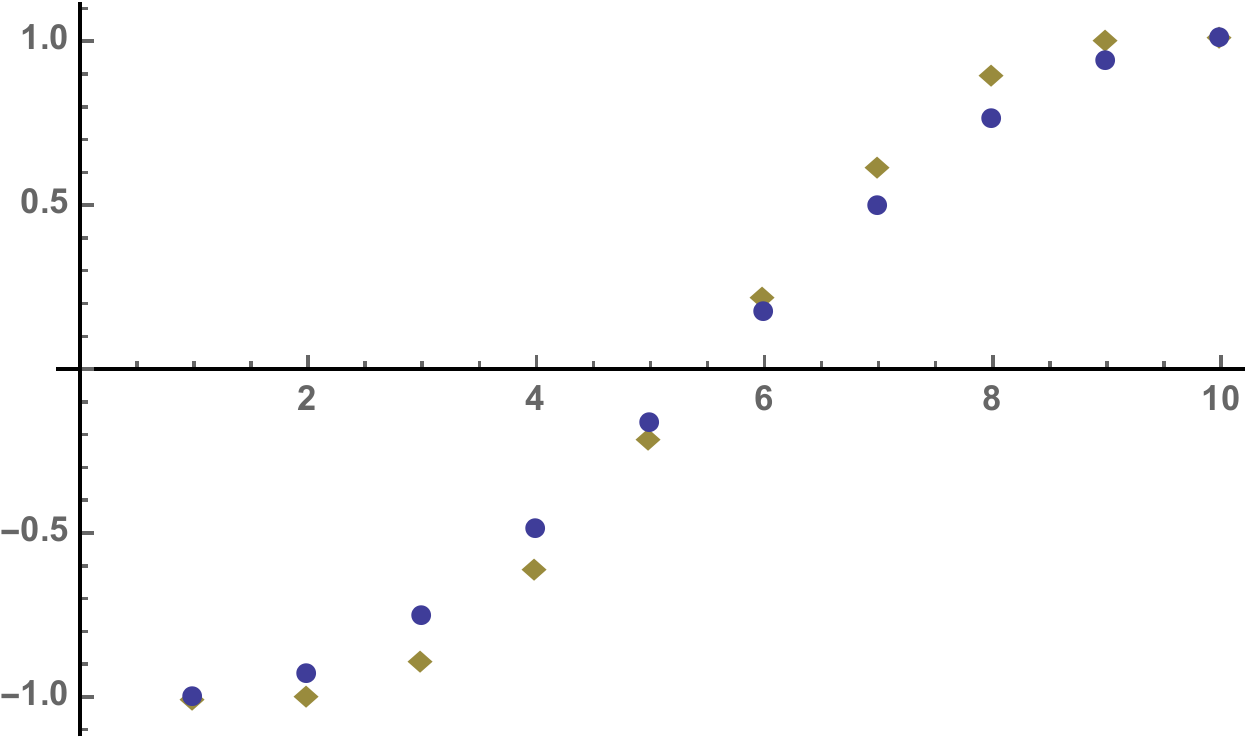}
              \caption{$n=10, \lambda=-3/4, m=10$. }
              \label{Ex6Zeros10}
          \end{figure}
      \end{minipage}
      \hspace{0.05\linewidth}
      \begin{minipage}{0.45\linewidth}
          \begin{figure}[H]
              \includegraphics[width=\linewidth]{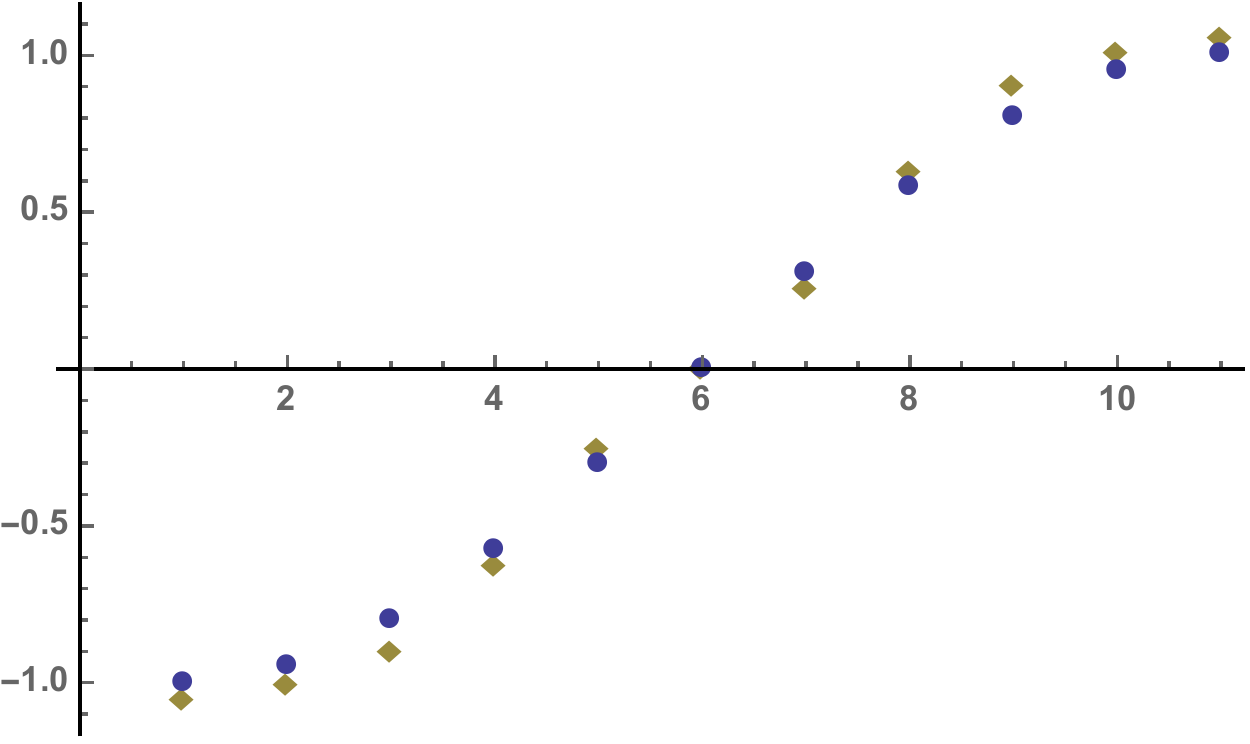}
              \caption{$n=10, \lambda=-3/4, m=11$.   }
              \label{Ex6Zeros11}
          \end{figure}
      \end{minipage}
  \end{minipage}

  \begin{minipage}{\linewidth}
      \centering
      \begin{minipage}{0.45\linewidth}
          \begin{figure}[H]
              \includegraphics[width=\linewidth]{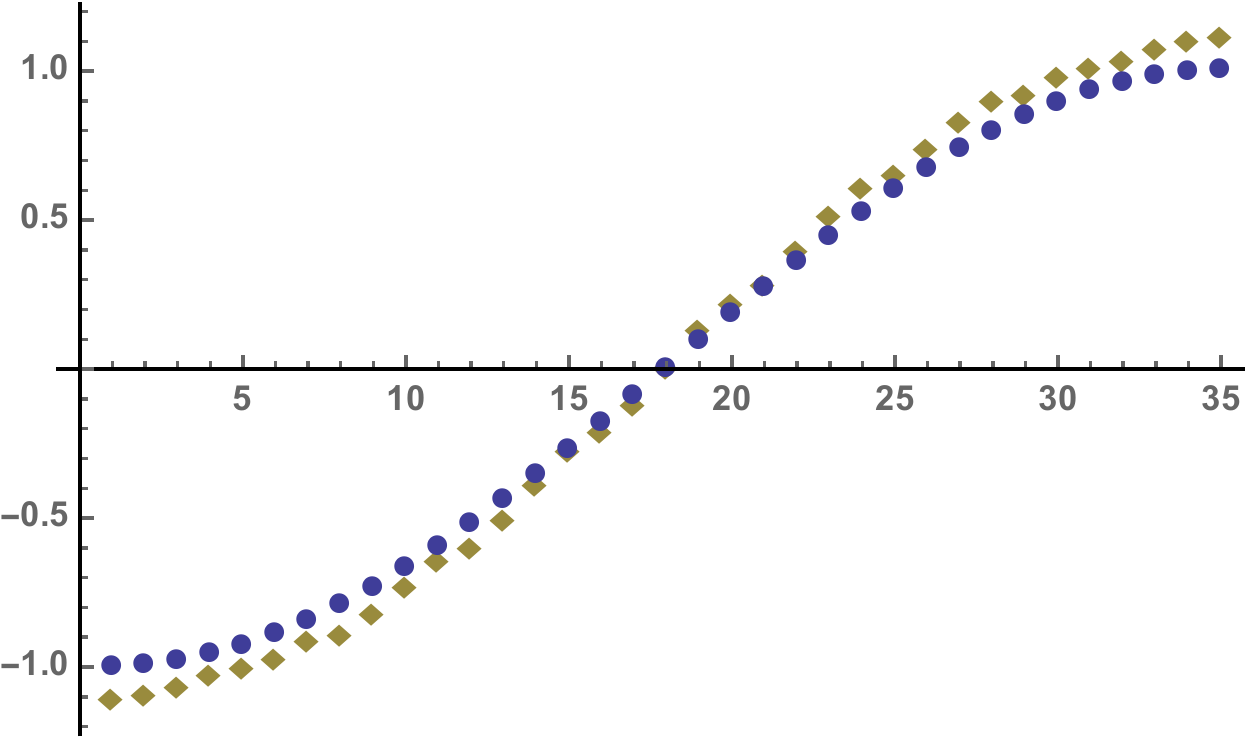}
              \caption{$n=10, \lambda=-3/4, m=35$.}
              \label{Ex6Zeros35}
          \end{figure}
      \end{minipage}
      \hspace{0.05\linewidth}
      \begin{minipage}{0.45\linewidth}
          \begin{figure}[H]
              \includegraphics[width=\linewidth]{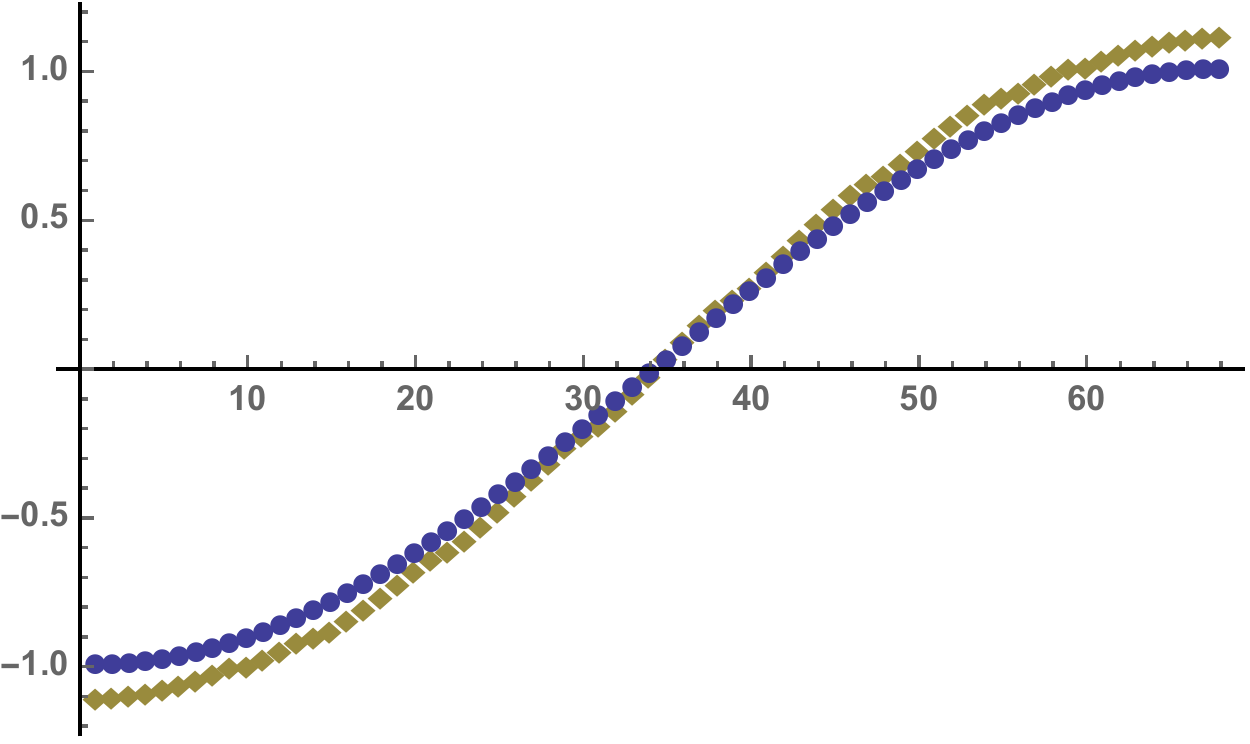}
              \caption{$n=10, \lambda=-3/4, m=68$.  }
              \label{Ex6Zeros68}
          \end{figure}
      \end{minipage}
  \end{minipage}

\smallskip

  \noindent \textbf{Example 5.5.}
Let $n=10$, $k=58$, $\sigma=2.$  Choose $\lambda=-1/4$  and  $a=\frac{4(2+\lambda)}{3(3+2\lambda)}=\frac{14}{15}$.
In Figures~\ref{Ex7Zeros8} through~\ref{Ex7Zeros67} the $y$-coordinates of the plotted points are the zeros of $D_m^\lambda$  and $C_m^\lambda$ for several integer values of $m$  between $8$ and $68$.

\begin{minipage}{\linewidth}
      \centering
      \begin{minipage}{0.45\linewidth}
          \begin{figure}[H]
              \includegraphics[width=\linewidth]{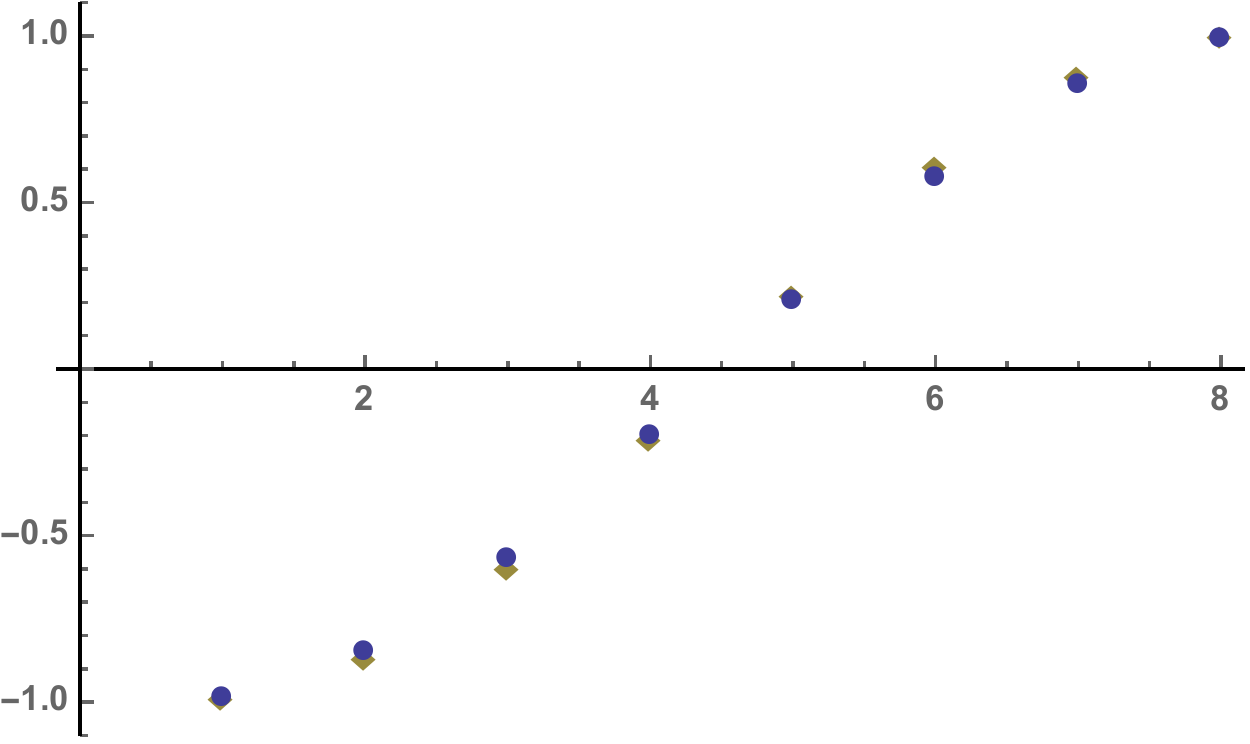}
              \caption{$n=10, \lambda=-1/4, m=8$.}
              \label{Ex7Zeros8}
          \end{figure}
      \end{minipage}
      \hspace{0.05\linewidth}
      \begin{minipage}{0.45\linewidth}
          \begin{figure}[H]
              \includegraphics[width=\linewidth]{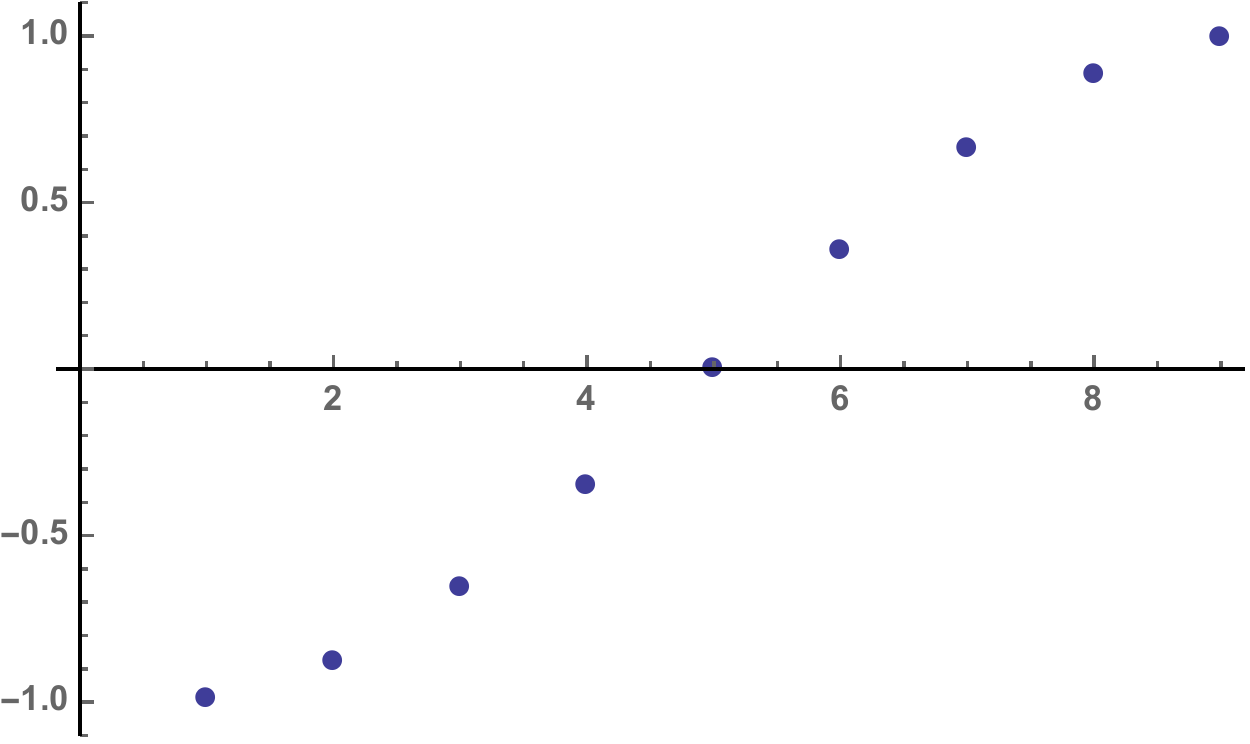}
              \caption{$n=10, \lambda=-1/4, m=9$. }
              \label{Ex7Zeros9}
          \end{figure}
      \end{minipage}
  \end{minipage}

  \begin{minipage}{\linewidth}
      \centering
      \begin{minipage}{0.45\linewidth}
          \begin{figure}[H]
              \includegraphics[width=\linewidth]{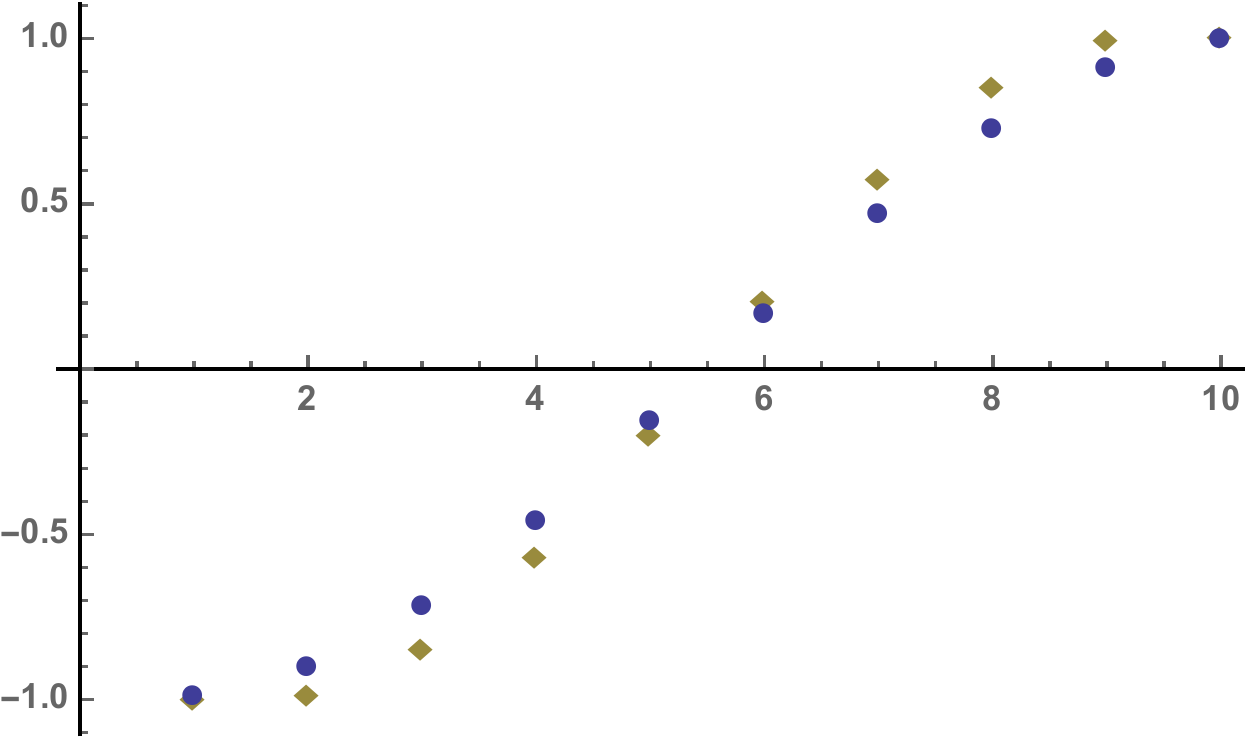}
              \caption{$n=10, \lambda=-1/4, m=10$.}
              \label{Ex7Zeros10}
          \end{figure}
      \end{minipage}
      \hspace{0.05\linewidth}
      \begin{minipage}{0.45\linewidth}
          \begin{figure}[H]
              \includegraphics[width=\linewidth]{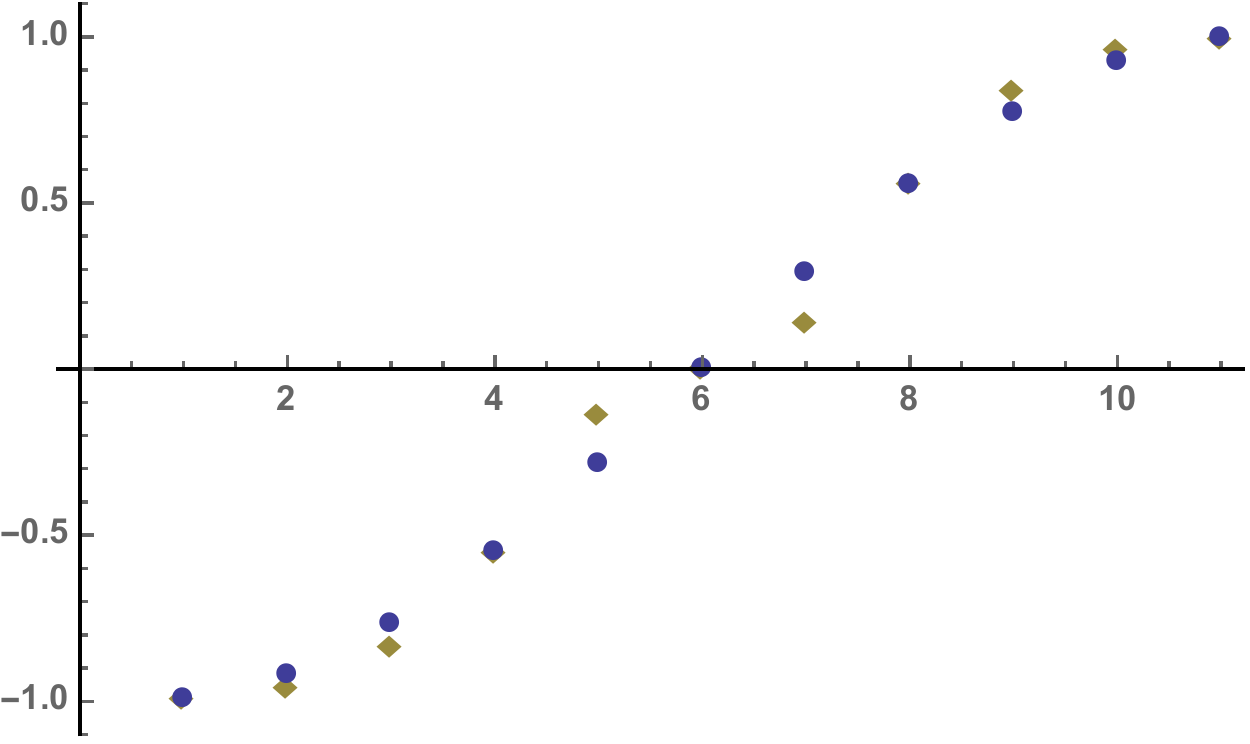}
              \caption{$n=10, \lambda=-1/4, m=11$.  }
              \label{Ex7Zeros11}
          \end{figure}
      \end{minipage}
  \end{minipage}

  \begin{minipage}{\linewidth}
      \centering
      \begin{minipage}{0.45\linewidth}
          \begin{figure}[H]
              \includegraphics[width=\linewidth]{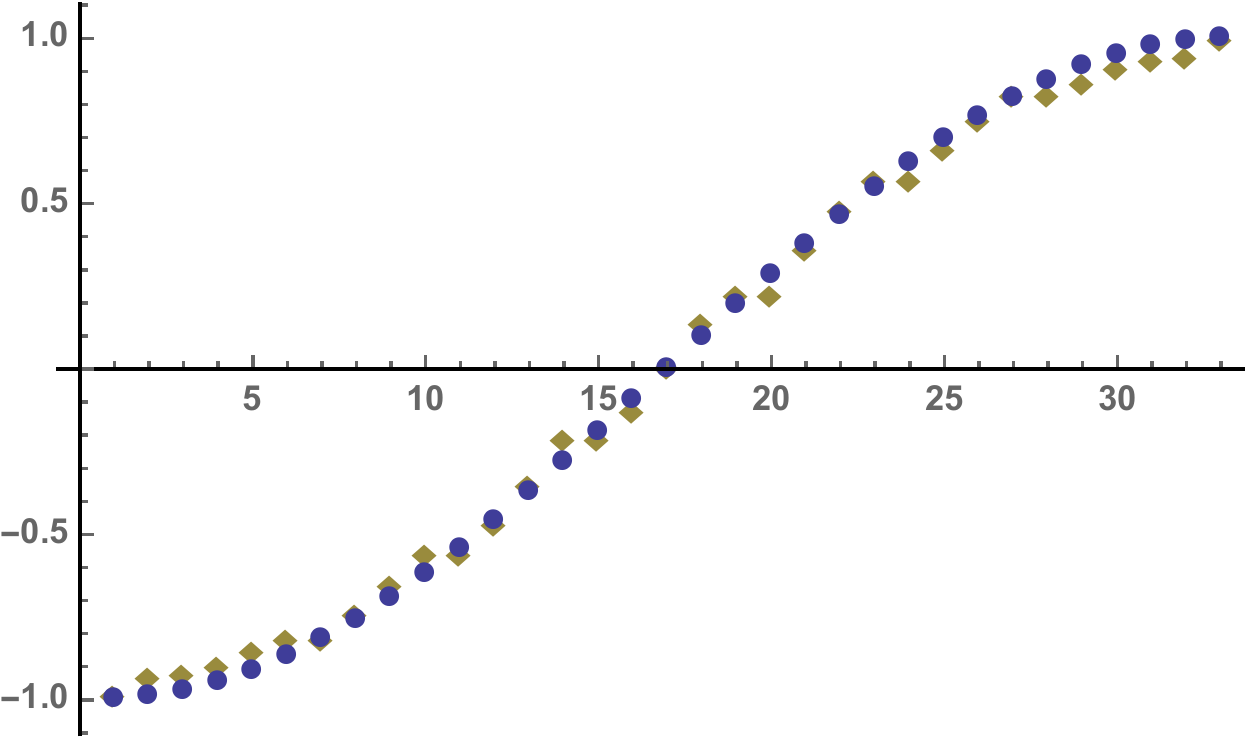}
              \caption{$n=10, \lambda=-1/4, m=34$.}
              \label{Ex7Zeros34}
          \end{figure}
      \end{minipage}
      \hspace{0.05\linewidth}
      \begin{minipage}{0.45\linewidth}
          \begin{figure}[H]
              \includegraphics[width=\linewidth]{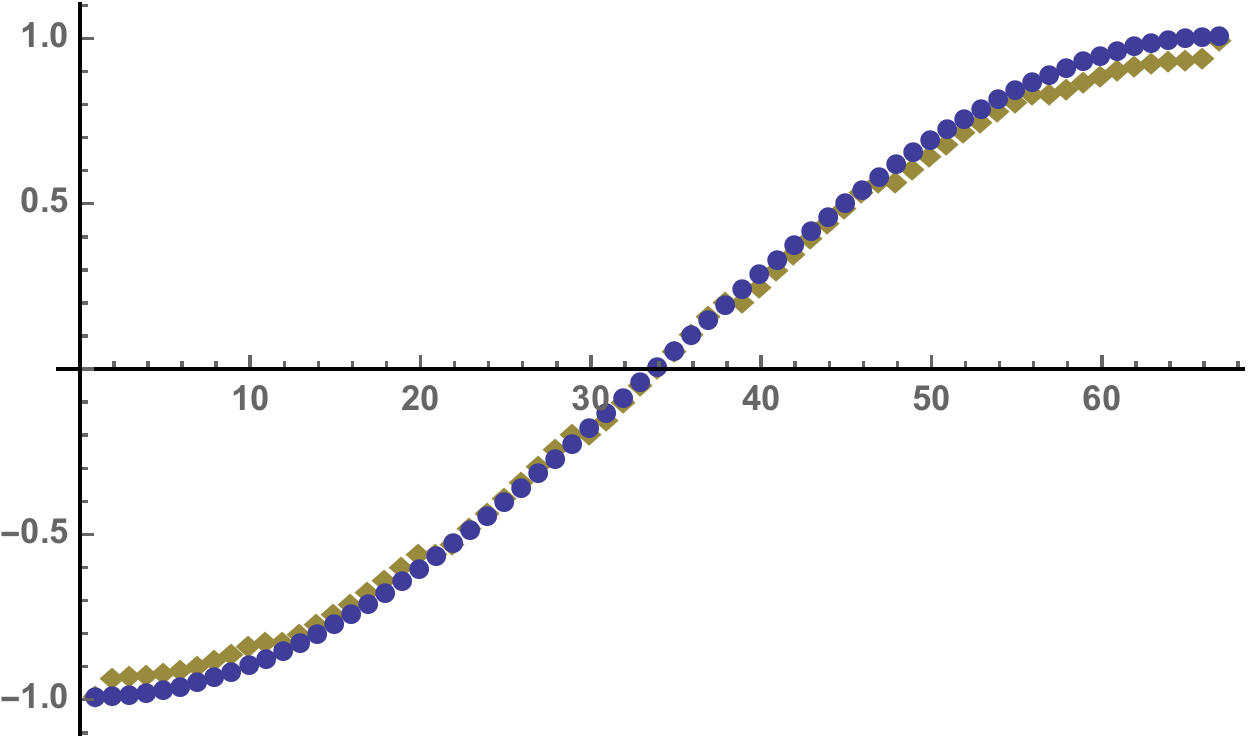}
              \caption{$n=10, \lambda=-1/4, m=67$. }
              \label{Ex7Zeros67}
          \end{figure}
      \end{minipage}
  \end{minipage}
  \smallskip

\noindent \textbf{Remark 5.2.} As $m$ increases, the zeros of $D_m^{\lambda}(x)$ and $C_m^{\lambda} (x)$ appear to be asymptotically equal for $\lambda = -1/4$. This is not unexpected since $\lambda = -1/4$ lies in the orthogonal range $\lambda > -1/2$ for ultraspherical polynomials. Note that the interval of orthogonality is $(-a,a)$,  where $a <1.$

\smallskip

    \noindent \textbf{Example 5.6.}
Let $n=5$, $k=18$, $\sigma=2.$  In Figures~\ref{Ex1LambdaM11_8} through~\ref{Ex1Lambda5_8}, the $y$-coordinates of the plotted points are the zeros of $D_m^\lambda$  and $C_m^\lambda$ for a selection of values of $\lambda \in (-\frac{3}{2}, +\infty)$ with  $\lambda \neq -1, 0, (2k-1)/2$, $k=0,1,2,\ldots$ where $m=23$ is fixed. We choose $a=\frac{4(2+\lambda)}{3(3+2\lambda)}$ if $\lambda<-1/2$ and $a=1$ if $\lambda>-1/2$; the zeros of  $D_m^\lambda$ and  $C_m^\lambda$ are contained in $(-a,a)$.

\noindent For  $ -3/2 < \lambda <  -1, $ the curves to which the zeros of  $D_m^{\lambda} (x)$ and  $C_m^{\lambda} (x)$ can be fitted are substantially different for some values of $m.$ As $\lambda$ approaches $-1/2$ from the left, the two curves are very similar, and,  as $\lambda>$ increases further,  the curves are almost identical.

\begin{minipage}{\linewidth}
      \centering
      \begin{minipage}{0.45\linewidth}
          \begin{figure}[H]
              \includegraphics[width=\linewidth]{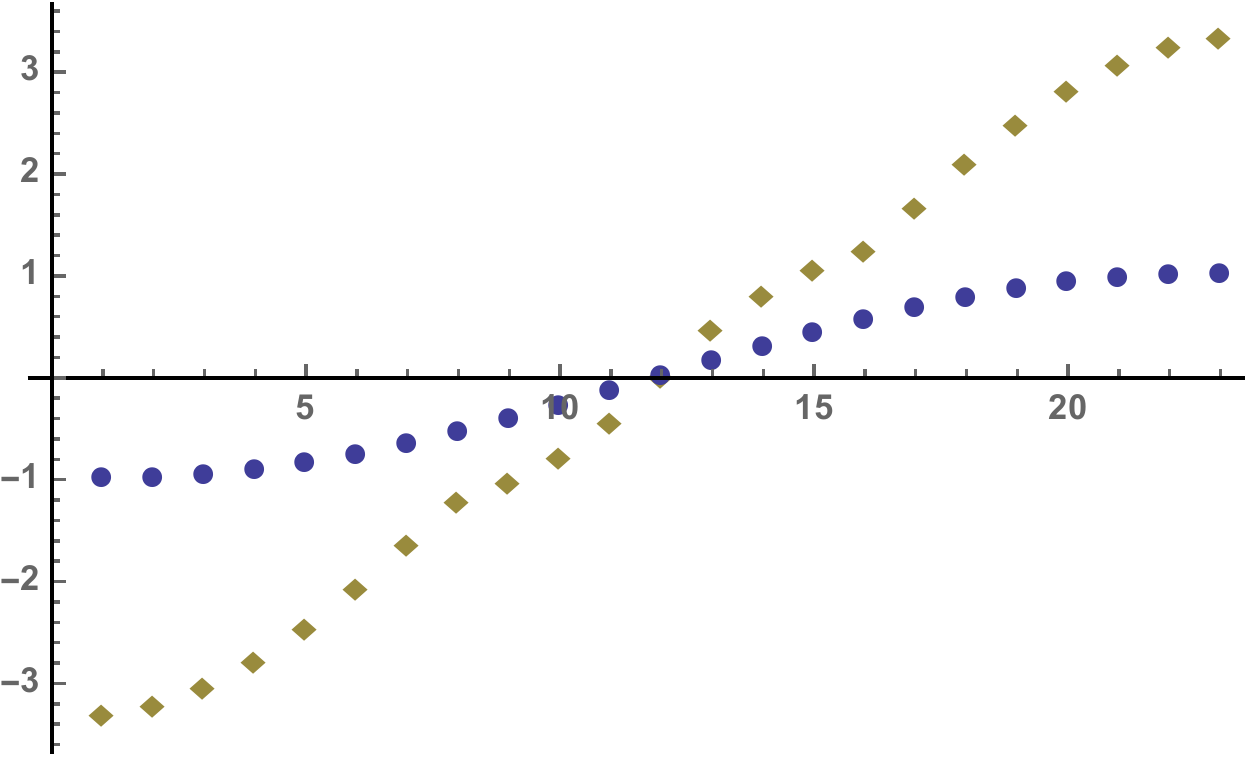}
              \caption{ $n=5, m=23, \lambda=-11/8, a=10/3$. }
              \label{Ex1LambdaM11_8}
          \end{figure}
      \end{minipage}
      \hspace{0.05\linewidth}
      \begin{minipage}{0.45\linewidth}
          \begin{figure}[H]
              \includegraphics[width=\linewidth]{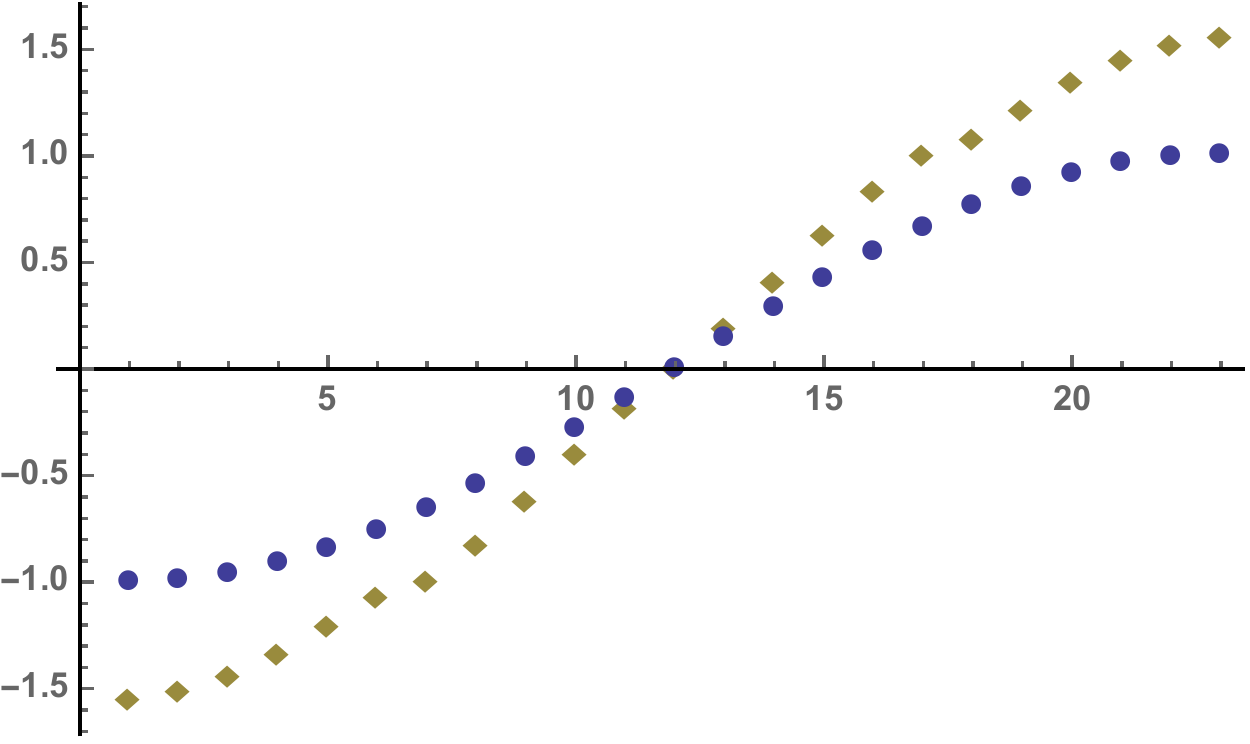}
              \caption{$n=5, m=23, \lambda=-9/8, a=14/9$.  }
              \label{Ex1LambdaM9_8}
          \end{figure}
      \end{minipage}
  \end{minipage}

  \begin{minipage}{\linewidth}
      \centering
      \begin{minipage}{0.45\linewidth}
          \begin{figure}[H]
              \includegraphics[width=\linewidth]{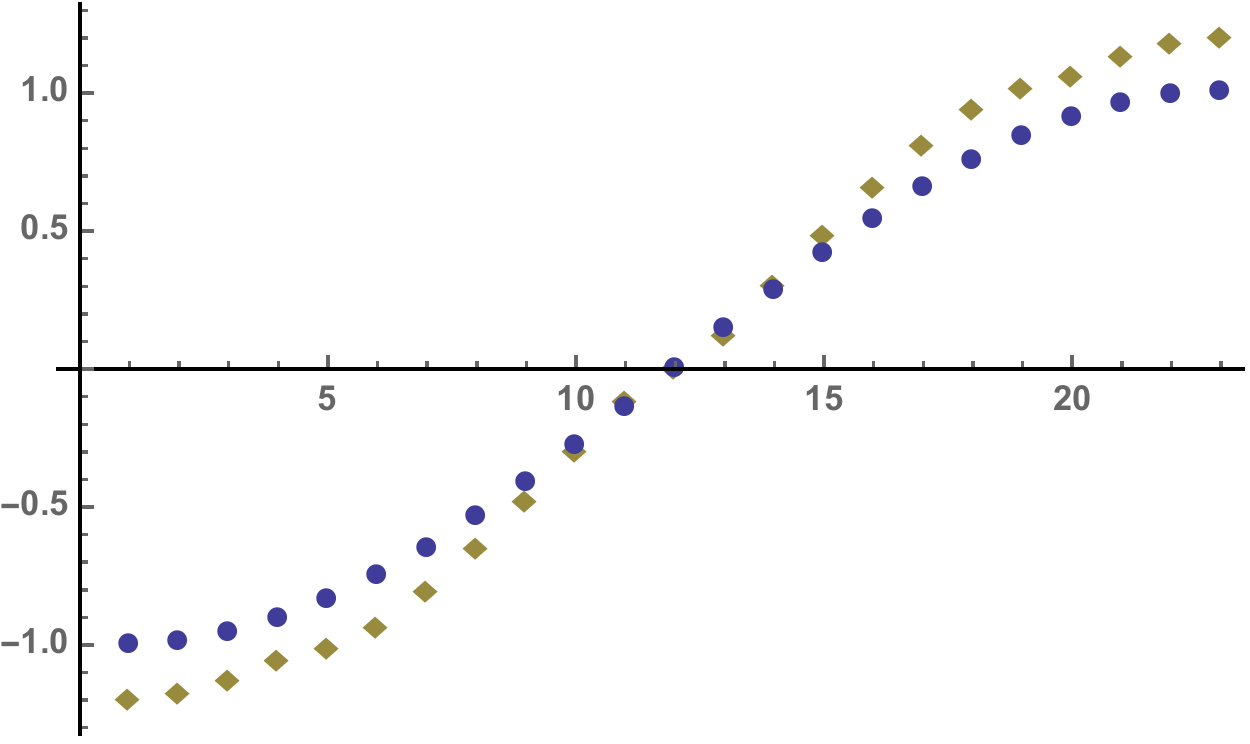}
              \caption{$n=5, m=23, \lambda=-7/8, a=6/5$. }
              \label{Ex1LambdaM7_8}
          \end{figure}
      \end{minipage}
      \hspace{0.05\linewidth}
      \begin{minipage}{0.45\linewidth}
          \begin{figure}[H]
              \includegraphics[width=\linewidth]{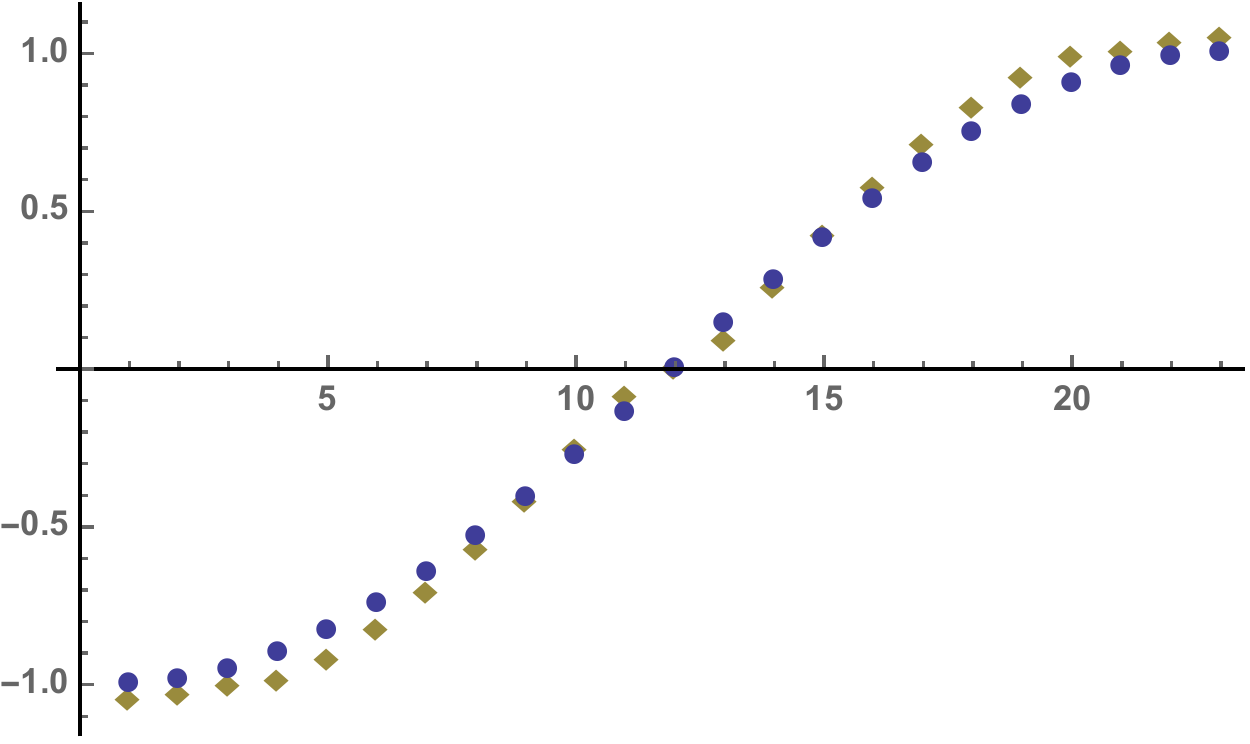}
              \caption{$n=5, m=23, \lambda=-5/8, a=22/21$. }
              \label{Ex1LambdaM5_8}
          \end{figure}
      \end{minipage}
  \end{minipage}

   \begin{minipage}{\linewidth}
      \centering
      \begin{minipage}{0.45\linewidth}
          \begin{figure}[H]
              \includegraphics[width=\linewidth]{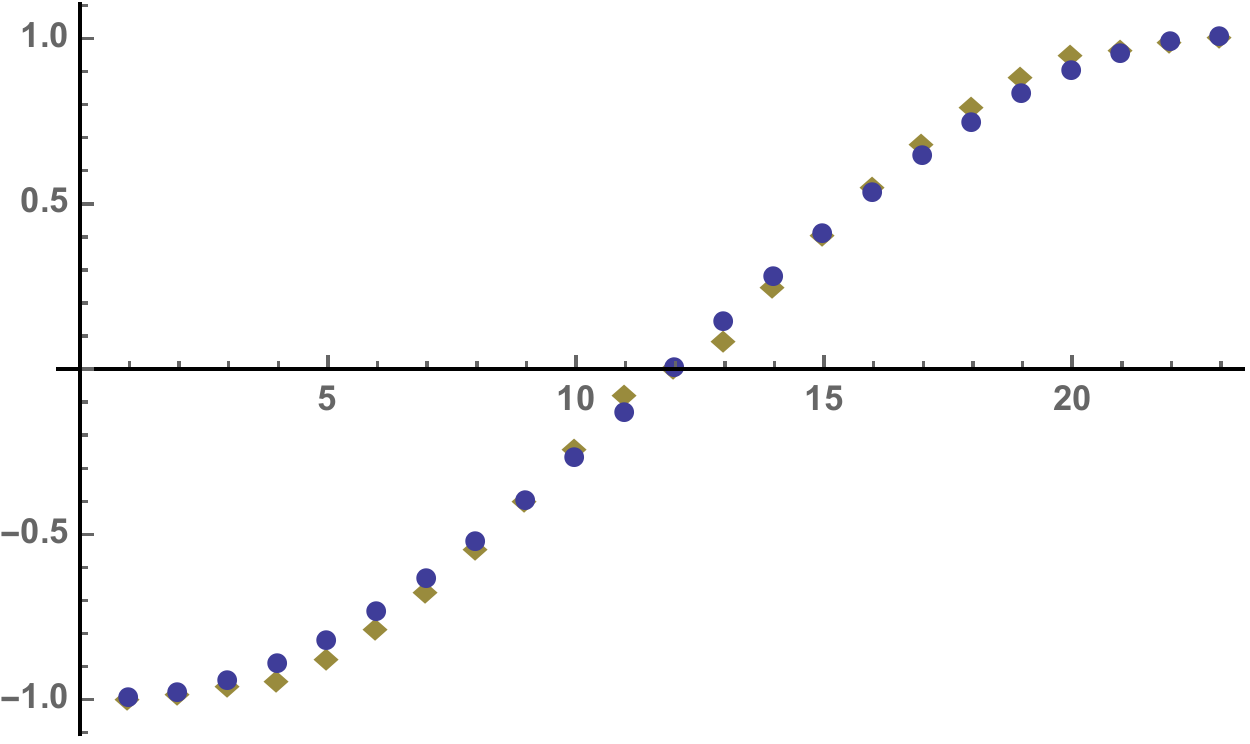}
              \caption{$n=5, m=23, \lambda=-3/8, a=1$. }
              \label{Ex1LambdaM3_8}
          \end{figure}
      \end{minipage}
      \hspace{0.05\linewidth}
      \begin{minipage}{0.45\linewidth}
          \begin{figure}[H]
              \includegraphics[width=\linewidth]{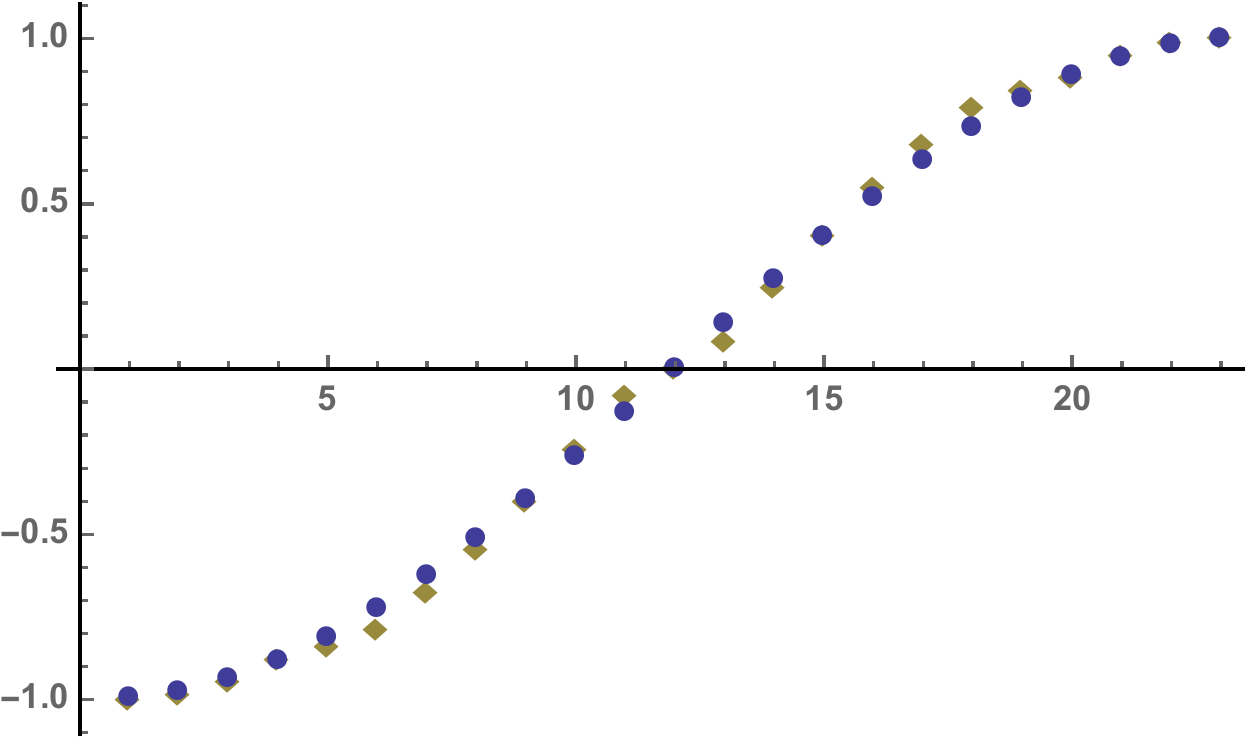}
              \caption{$n=5, m=23, \lambda=1/8, a=1$. }
              \label{Ex1Lambda1_8}
          \end{figure}
      \end{minipage}
  \end{minipage}

   \begin{minipage}{\linewidth}
      \centering
      \begin{minipage}{0.45\linewidth}
          \begin{figure}[H]
              \includegraphics[width=\linewidth]{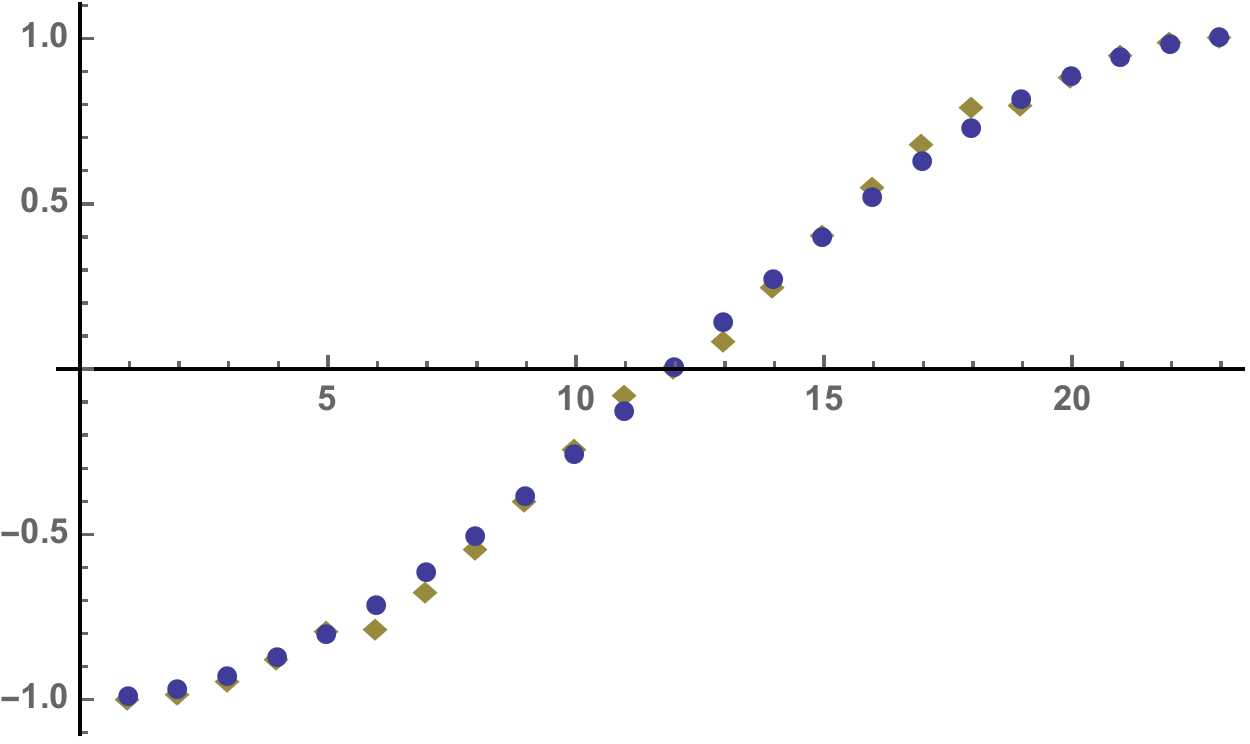}
              \caption{$n=5, m=23, \lambda=3/8, a=1$.}
              \label{Ex1Lambda3_8}
          \end{figure}
      \end{minipage}
      \hspace{0.05\linewidth}
      \begin{minipage}{0.45\linewidth}
          \begin{figure}[H]
              \includegraphics[width=\linewidth]{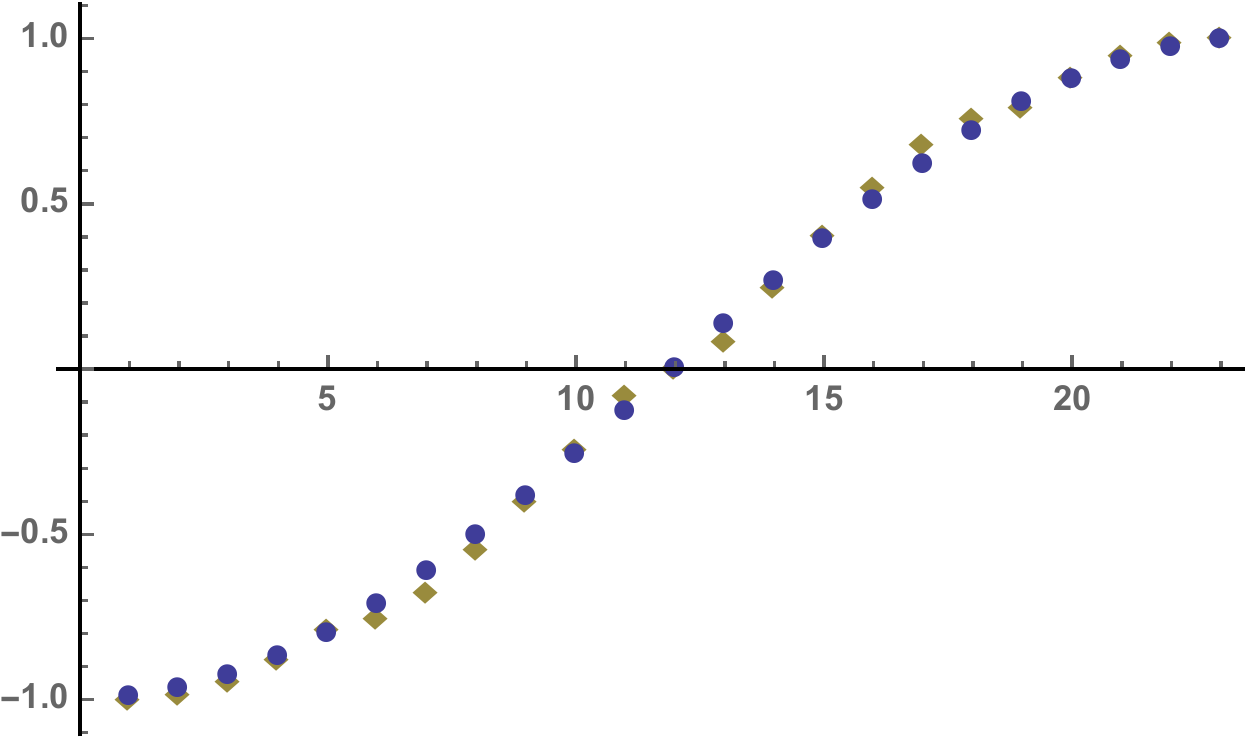}
              \caption{$n=5, m=23, \lambda=5/8, a=1$. }
              \label{Ex1Lambda5_8}
          \end{figure}
      \end{minipage}
  \end{minipage}

\medskip

\section{\large{Acknowledgements}}

\noindent Kathy Driver would like to express her thanks to the Mathematics Department, University of Colorado, Colorado Springs, for its hospitality during her visit in Spring 2018, during which the work on this paper began.

 \end{document}